%
\magnification=\magstep1   
\input amstex
\UseAMSsymbols
\input pictex

\vsize=23truecm
\NoBlackBoxes
\parindent=18pt
  
   \font\rmk=cmr8      \font\ttk=cmtt8

\font\gross=cmbx10 scaled\magstep1 

\def\odd{\text{odd}}
\def\ev{\text{ev}}

\def\arr#1#2{\arrow <1.5mm> [0.25,0.75] from #1 to #2}

   

 \def\Rahmenbio#1%
   {$$\vbox{\hrule\hbox%
                  {\vrule%
                       \hskip0.5cm%
                            \vbox{\vskip0.3cm\relax%
                               \hbox{$\displaystyle{#1}$}%
                                  \vskip0.3cm}%
                       \hskip0.5cm%
                  \vrule}%
           \hrule}$$}        

\def\Rahmen#1%
   {\centerline{\vbox{\hrule\hbox%
                  {\vrule%
                       \hskip0.5cm%
                            \vbox{\vskip0.3cm\relax%
                               \hbox{{#1}}%
                                  \vskip0.3cm}%
                       \hskip0.5cm%
                  \vrule}%
           \hrule}}}

\centerline{\gross The Fibonacci partition triangles}
	\medskip
\centerline{Philipp Fahr, Claus Michael Ringel}
	\medskip

	\bigskip

{\narrower\narrower{\bf Abstract.}
In two previous papers we have presented partition formulae for the 
Fibonacci numbers motivated by the appearance of the Fibonacci numbers
in the representation theory of the 
$3$-Kronecker quiver and its universal cover, the $3$-regular tree.
Here we show that the basic information can be rearranged in 
two triangles. They are quite similar to the Pascal
triangle of the binomial coefficients, but in contrast to the additivity
rule for the Pascal triangle, we now deal with additivity along ``hooks'', or, 
equivalently, with additive functions for valued translation quivers.
As for the Pascal triangle, 
we see that the numbers in these Fibonacci partition triangles 
are given by evaluating polynomials. We show that
the two triangles can be obtained from each other by looking at 
differences of numbers, it is sufficient to take differences
along arrows and knight's moves.\par}
	\bigskip
The aim of the paper is to rearrange the positive integers which
are used in the partition formulae for the Fibonacci numbers as considered
in [FR1,FR2].
For the even-index Fibonacci numbers we obtain a proper triangle which
we call the even-index Fibonacci partition triangle. Second, what we call
the odd-index Fibonacci partition triangle actually is a triangle only after
removing one number (but it seems worthwhile to take this additional 
position into account). These arrangements of integers are quite similar to the
Pascal triangle of the binomial coefficients. In particular, we will show
that the diagonals are given by evaluating polynomials. 

Let us recall that our work on the partition formulae was based on
the appearance of the Fibonacci numbers in the representation
theory of quivers and the present paper again relies on concepts 
which have been developed in this context.
Namely, both triangles turn out to show additive functions
on valued translation quivers. Translation quivers have been
considered frequently in the representation theory of quivers and 
finite-dimensional algebras.
The valued translation quivers can be used in order to describe 
module categories, but those with non-trivial valuation (as in
the case of the Fibonacci partition triangles) have  seldom be seen to be of
importance when dealing with quivers.

In sections 1 and 2 we will exhibit the even-index, and the odd-index triangle,
respectively. 
Our main task will be to show in which 
way the two triangles can be obtained from each other, see sections 3 and 4.
The relationship which we will encounter shows that the two
triangles are intimately connected. The proof provided here
relies on the categorification of the Fibonacci pairs given in [FR2], it will
be given in section 5. The final section 6 provides some further remarks
and open questions.

	\bigskip
The investigation is based on the Fibonacci partition formulae established
in [FR1] for the Fibonacci numbers with even index, 
and in the  PhD thesis [F] of Fahr for those with odd index, 
see also [FR2], but also
on further discussions of the authors during the time when
Fahr was a PhD student at Bielefeld. The final version was
written by Ringel.

\vfill\eject
{\bf 1. The even-index Fibonacci partition triangle}

$$
{\beginpicture
\setcoordinatesystem units <.8cm,1cm>
\multiput{$\ssize (2,1)$} at -1.5  2  
 -2.5  2  -3.5  2  -4.5  2  -5.5  2  -6.5  2  -7.5  2  -8.5  2  -9.5  2  -10.5  2 /

\put{$\ssize (3,1)$} at  -.5  2 

\multiput{$1$} at 0 0  -1 -1 -2 -2  -3 -3  -4 -4  -5 -5  -6 -6  -7 -7  -8 -8  -9 -9 
 -10 -10 /
\multiput{$\cdots$} at -11.5 2  -6 1 /

\put{Valuation:} [l] at -14.4 2.05

\setdots <1mm>
\plot 0 1.5  0 .5 /
\plot -1 2.5  -1 .5 /
\plot -2 2.5  -2 .5 /
\plot -3 2.5  -3 .5 /

\setsolid

\put{$2$} at 0 -2 
\put{$3$} at -1 -3 
\put{$4$} at -2 -4 
\put{$5$} at -3 -5 
\put{$6$} at -4 -6 
\put{$7$} at -5 -7 
\put{$8$} at -6 -8 
\put{$9$} at -7 -9 
\put{$10$} at -8 -10 

\put{$7$} at 0 -4
\put{$12$} at -1 -5
\put{$18$} at -2 -6 
\put{$25$} at -3 -7 
\put{$33$} at -4 -8 
\put{$42$} at -5 -9 
\put{$52$} at -6 -10 

\put{$29$} at 0 -6
\put{$53$} at -1 -7 
\put{$85$} at -2 -8 
\put{$126$} at -3 -9 
\put{$177$} at -4 -10 

\put{$130$} at 0 -8
\put{$247$} at -1 -9 
\put{$414$} at -2 -10

\put{$611$} at 0 -10

\put{$1$} at -11 -11
\put{$11$} at -9 -11
\put{$63$} at -7 -11
\put{$239$} at -5 -11
\put{$642$} at -3 -11
\put{$1192$} at -1 -11
\put{$1$} at -12 -12
\put{$12$} at -10 -12
\put{$75$} at -8 -12
\put{$313$} at -6 -12
\put{$943$} at -4 -12
\put{$2062$} at -2 -12
\put{$2965$} at 0 -12

\arr{-.3 -.3}{-.7 -.7}

\arr{-1.3 -1.3}{-1.7 -1.7}
\arr{-.7 -1.3}{-.3 -1.7}

\arr{-2.3 -2.3}{-2.7 -2.7}
\arr{-1.7 -2.3}{-1.3 -2.7}
\arr{-.3 -2.3}{-.7 -2.7}

\arr{-3.3 -3.3}{-3.7 -3.7}
\arr{-2.7 -3.3}{-2.3 -3.7}
\arr{-1.3 -3.3}{-1.7 -3.7}
\arr{-.7 -3.3}{-.3 -3.7}

\arr{-4.3 -4.3}{-4.7 -4.7}
\arr{-3.7 -4.3}{-3.3 -4.7}
\arr{-2.3 -4.3}{-2.7 -4.7}
\arr{-1.7 -4.3}{-1.3 -4.7}
\arr{-.3 -4.3}{-.7 -4.7}

\arr{-5.3 -5.3}{-5.7 -5.7}
\arr{-4.7 -5.3}{-4.3 -5.7}
\arr{-3.3 -5.3}{-3.7 -5.7}
\arr{-2.7 -5.3}{-2.3 -5.7}
\arr{-1.3 -5.3}{-1.7 -5.7}
\arr{-.7 -5.3}{-.3 -5.7}

\arr{-6.3 -6.3}{-6.7 -6.7}
\arr{-5.7 -6.3}{-5.3 -6.7}
\arr{-4.3 -6.3}{-4.7 -6.7}
\arr{-3.7 -6.3}{-3.3 -6.7}
\arr{-2.3 -6.3}{-2.7 -6.7}
\arr{-1.7 -6.3}{-1.3 -6.7}
\arr{-.3 -6.3}{-.7 -6.7}

\arr{-7.3 -7.3}{-7.7 -7.7}
\arr{-6.7 -7.3}{-6.3 -7.7}
\arr{-5.3 -7.3}{-5.7 -7.7}
\arr{-4.7 -7.3}{-4.3 -7.7}
\arr{-3.3 -7.3}{-3.7 -7.7}
\arr{-2.7 -7.3}{-2.3 -7.7}
\arr{-1.3 -7.3}{-1.7 -7.7}
\arr{-.7 -7.3}{-.3 -7.7}

\arr{-8.3 -8.3}{-8.7 -8.7}
\arr{-7.7 -8.3}{-7.3 -8.7}
\arr{-6.3 -8.3}{-6.7 -8.7}
\arr{-5.7 -8.3}{-5.3 -8.7}
\arr{-4.3 -8.3}{-4.7 -8.7}
\arr{-3.7 -8.3}{-3.3 -8.7}
\arr{-2.3 -8.3}{-2.7 -8.7}
\arr{-1.7 -8.3}{-1.3 -8.7}
\arr{-.3 -8.3}{-.7 -8.7}

\arr{-9.3 -9.3}{-9.7 -9.7}
\arr{-8.7 -9.3}{-8.3 -9.7}
\arr{-7.3 -9.3}{-7.7 -9.7}
\arr{-6.7 -9.3}{-6.3 -9.7}
\arr{-5.3 -9.3}{-5.7 -9.7}
\arr{-4.7 -9.3}{-4.3 -9.7}
\arr{-3.3 -9.3}{-3.7 -9.7}
\arr{-2.7 -9.3}{-2.3 -9.7}
\arr{-1.3 -9.3}{-1.7 -9.7}
\arr{-.7 -9.3}{-.3 -9.7}

\arr{-10.3 -10.3}{-10.7 -10.7}
\arr{-9.7 -10.3}{-9.3 -10.7}
\arr{-8.3 -10.3}{-8.7 -10.7}
\arr{-7.7 -10.3}{-7.3 -10.7}
\arr{-6.3 -10.3}{-6.7 -10.7}
\arr{-5.7 -10.3}{-5.3 -10.7}
\arr{-4.3 -10.3}{-4.7 -10.7}
\arr{-3.7 -10.3}{-3.3 -10.7}
\arr{-2.3 -10.3}{-2.7 -10.7}
\arr{-1.7 -10.3}{-1.3 -10.7}
\arr{-.3 -10.3}{-.7 -10.7}

\arr{-11.3 -11.3}{-11.7 -11.7}
\arr{-10.7 -11.3}{-10.3 -11.7}
\arr{-9.3 -11.3}{-9.7 -11.7}
\arr{-8.7 -11.3}{-8.3 -11.7}
\arr{-7.3 -11.3}{-7.7 -11.7}
\arr{-6.7 -11.3}{-6.3 -11.7}
\arr{-5.3 -11.3}{-5.7 -11.7}
\arr{-4.7 -11.3}{-4.3 -11.7}
\arr{-3.3 -11.3}{-3.7 -11.7}
\arr{-2.7 -11.3}{-2.3 -11.7}
\arr{-1.3 -11.3}{-1.7 -11.7}
\arr{-.7 -11.3}{-.3 -11.7}

\plot -12.3 -12.3 -12.7 -12.7 /
\plot -11.7 -12.3 -11.3 -12.7 /
\plot -10.3 -12.3 -10.7 -12.7 /
\plot -9.7 -12.3 -9.3 -12.7 /
\plot -8.3 -12.3 -8.7 -12.7 /
\plot -7.7 -12.3 -7.3 -12.7 /
\plot -6.3 -12.3 -6.7 -12.7 /
\plot -5.7 -12.3 -5.3 -12.7 /
\plot -4.3 -12.3 -4.7 -12.7 /
\plot -3.7 -12.3 -3.3 -12.7 /
\plot -2.3 -12.3 -2.7 -12.7 /
\plot -1.7 -12.3 -1.3 -12.7 /
\plot  -.3 -12.3  -.7 -12.7 /

\linethickness.7mm
\putrule from 0 -.3 to 0 -1.7
\putrule from 0 -2.3 to 0 -3.7
\putrule from 0 -4.3 to 0 -5.7
\putrule from 0 -6.3 to 0 -7.7
\putrule from 0 -8.3 to 0 -9.7
\putrule from 0 -10.3 to 0 -11.7
\putrule from 0 -12.3 to 0 -12.7

\putrule from 0 2.3 to 0 1.6

\put{$\ssize 0$} at 2.5 0
\put{$\ssize 1$} at 2.5 -1
\put{$\ssize 2$} at 2.5 -2
\put{$\ssize 3$} at 2.5 -3
\put{$\ssize 4$} at 2.5 -4
\put{$\ssize 5$} at 2.5 -5
\put{$\ssize 6$} at 2.5 -6
\put{$\ssize 7$} at 2.5 -7
\put{$\ssize 8$} at 2.5 -8
\put{$\ssize 9$} at 2.5 -9
\put{$\ssize 10$} at 2.5 -10
\put{$\ssize 11$} at 2.5 -11
\put{$\ssize 12$} at 2.5 -12

\put{$\ssize f_2$} at 1.5 0
\put{$\ssize f_4$} at 1.5 -1
\put{$\ssize f_6$} at 1.5 -2
\put{$\ssize f_8$} at 1.5 -3
\put{$\ssize f_{10}$} at 1.5 -4
\put{$\ssize f_{12}$} at 1.5 -5
\put{$\ssize f_{14}$} at 1.5 -6
\put{$\ssize f_{16}$} at 1.5 -7
\put{$\ssize f_{18}$} at 1.5 -8
\put{$\ssize f_{20}$} at 1.5 -9
\put{$\ssize f_{22}$} at 1.5 -10
\put{$\ssize f_{24}$} at 1.5 -11
\put{$\ssize f_{26}$} at 1.5 -12

\put{$d_0(t)$} at -14 -14
\put{$d_1(t)$} at -12 -14
\put{$d_2(t)$} at -10 -14
\put{$d_3(t)$} at -8 -14
\put{$d_4(t)$} at -6 -14
\put{$d_5(t)$} at -4 -14
\put{$d_6(t)$} at -2 -14
\put{$t$} at 2.5 -14
\put{$f_{2t+2}$} at 1.5 -14

\arr{-13.3 -13.3}{-13.7 -13.7}
\arr{-11.3 -13.3}{-11.7 -13.7}
\arr{-9.3 -13.3}{-9.7 -13.7}
\arr{-7.3 -13.3}{-7.7 -13.7}
\arr{-5.3 -13.3}{-5.7 -13.7}
\arr{-3.3 -13.3}{-3.7 -13.7}
\arr{-1.3 -13.3}{-1.7 -13.7}

\setdots <1mm>
\plot -14.5 0  -3 0 /
\plot -14.5 -1  -4 -1 /
\plot -14.5 -2  -5 -2 /
\plot -14.5 -3  -6 -3 /
\plot -14.5 -4  -7 -4 /
\plot -14.5 -5  -8 -5 /
\plot -14.5 -6  -9 -6 /
\plot -14.5 -7  -10 -7 /
\plot -14.5 -8  -11 -8 /
\plot -14.5 -9  -12 -9 /
\plot -14.5 -10  -13 -10 /

\put{$P_1$} at -14 -.5
\put{$P_2$} at -14 -1.5
\put{$P_3$} at -14 -2.5
\put{$P_4$} at -14 -3.5
\put{$P_5$} at -14 -4.5
\put{$P_6$} at -14 -5.5
\put{$P_7$} at -14 -6.5
\put{$P_8$} at -14 -7.5
\put{$P_9$} at -14 -8.5
\put{$P_{10}$} at -14 -9.5

\endpicture}
$$
\vfill\eject

\noindent
First, let us explain the numbers displayed in the triangle as well as the two
columns on the right. The rows of the triangle are 
indexed by $t = 0,1,2,\dots$ as shown on the right.
In any row, say the row $t$, the entries will be labeled, from left to right, by
$d_0(t), d_1(t), \dots, d_i(t), \dots$, with $i \le \frac t2;$ for example,
for $t = 5$ we have 
$$
 d_0(5) = 1,\quad d_1(5) = 5, \quad d_2(5) = 12.
$$
For some calculations it seems convenient to define $d_i(t)$ 
for all $i\in \Bbb Z$ as follows: we let $d_i(t) = 0$ for $i < 0$ and we
define $d_i(t) = d_{t-i}(t)$ for $i > \frac t2$ (thus $d_i(t) = d_{t-i}(t)$
for all $i\in \Bbb Z$).

The entries $d_i(t)$ displayed in the triangle 
are calculated inductively as follows: We start with $d_0(t) = 1$ for all $t\ge 0$.
Now let $i \ge 1.$ Then, for $1 \le i < \frac t2$ (and $t\ge 3$) let
$$
 d_i(t) = 2d_{i-1}(t-1)+d_i(t-1) - d_{i-1}(t-2),
$$
whereas for $i = \frac t2$ (thus $t = 2i$)
$$
 d_i(2i) = 3d_{i-1}(2i-1)-d_{i-1}(2i-2).
$$
So using the convention that $d_i(t) = d_{t-i}(t)$ for $i> \frac t2$, we
see that the rule for $i = \frac t2$ is the same as that for $i < \frac t2$.
	\medskip

The entries $d_i(t)$ displayed in the triangle have been considered already in
the paper [FR1] (but labeled differently): they are derived from 
the numbers $a_s[j]$ of [FR1] according to the rule
$$
 d_i(t) = a_{\lceil t/2 \rceil}[t-2i]
$$
for $0 \le i \le \frac t2$
(given any real number $\alpha$, we denote by $\lceil \alpha \rceil$ 
the smallest integer $z$ with $\alpha \le z$).
	\medskip
Using the notation $d_i(t)$, the partition formulae of [FR1] can be written in
a unified way, as follows: 
	\medskip
{\bf Partition formula for the Fibonacci numbers $f_n$ with even index $n$:}
	\medskip
\Rahmen{$
 3\Biggl(\sum\limits_{0\le i < t/2} 2^{t-2i-1}d_i(t)\Biggr) +  d_{t/2}(t) \ =  f_{2t+2}
$}

\noindent
where we set $d_\alpha(t) = 0$ if $\alpha\notin \Bbb Z.$
This shows that {\it the row $t$ of the triangle yields the Fibonacci number $f_{2t+2}$}.
Thus, we have added the even-index Fibonacci numbers in a column on the right.
As examples, let us look at the rows $t= 5$ and $t=6$, they yield
$$
\align
 3(2^4\cdot 1 + 2^2\cdot 5 + 2^0\cdot 12) \qquad \ \ &= 144 =  f_{12}
      \tag{$t=5$}\cr
 3(2^5\cdot 1 + 2^3\cdot 6 + 2^1\cdot 18) + 29 \  &= 377 = f_{14} \tag{$t=6$}
\endalign
$$
	\bigskip

The numbers of the triangle are connected by arrows: {\bf we consider
the triangle} as a quiver, or better {\bf as a valued translation quiver.}
Let us recall the definition (see for example [HPR]).

A {\it quiver $\Gamma = (\Gamma_0,\Gamma_1)$ without multiple arrows} is given by 
a set $\Gamma_0$, called the set of {\it vertices,}
and a subset $\Gamma_1$ of $\Gamma_0\times \Gamma_0.$ An element $\alpha = (x,y)$ of $\Gamma_1$ 
with $x,y\in \Gamma_0$ is called an 
{\it arrow,} and usually one writes $\alpha\:x \to y$
and calls $x$ the {\it starting} vertex, $y$ the {\it terminal} vertex of $\alpha$. 
Given vertices $x, z$ of $\Gamma$, we denote by $x^+$ the set of vertices $y$ with
an arrow $x \to y$, and by $z^-$ the set of vertices $y$ with an arrow
$y \to z$. 

A {\it translation quiver} $\Gamma = (\Gamma_0,\Gamma_1,\tau)$ is
given by a quiver $(\Gamma_0,\Gamma_1)$ without multiple arrows, 
a subset $\Gamma_0^p \subseteq \Gamma_0$ and 
an injective function $\tau\:(\Gamma_0\setminus\Gamma_0^p) \to \Gamma_0$  
such that for any vertex $z$  in $\Gamma_0\setminus\Gamma_0^p$, 
one has $(\tau z)^+ = z^-$. The function $\tau$
is called the {\it translation,} the vertices in 
$\Gamma_0^p$ are said to be the {\it projective} vertices,
those not in the image of $\tau$ the {\it injective} vertices.
Given a non-projective vertex $z$, one says that $\tau z, z^-,z$ is the mesh
ending in $z$; in the examples which we consider, these sets $z^-$ consist of
either one or two elements, thus we deal with meshes of the following form,
where $x = \tau z$:
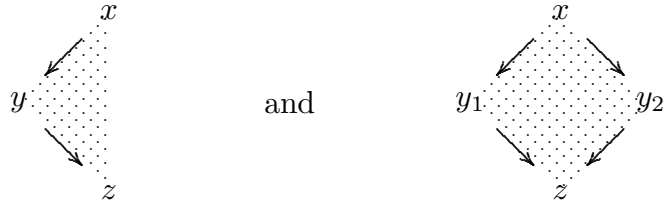
$$
{\beginpicture
\setcoordinatesystem units <1cm,1cm>
\put{\beginpicture
\put{$x$} at 1 2
\put{$y$} at 0 1
\put{$z$} at 1 0
\arr{0.7 1.7}{0.3 1.3}
\arr{0.3  .7}{0.7  .3}
\setshadegrid span <0.8mm>
\vshade 0 1 1 <z,z,z,z> 1 0 2  /
\endpicture} at 0 0 
\put{\beginpicture
\put{$x$} at 1 2
\put{$y_1$} at 0 1
\put{$y_2$} at 2 1
\put{$z$} at 1 0
\arr{0.7 1.7}{0.3 1.3}
\arr{1.3 1.7}{1.7 1.3}
\arr{0.3  .7}{0.7  .3}
\arr{1.7 .7}{1.3 .3}
\setshadegrid span <0.8mm>
\vshade 0 1 1 <z,z,z,z> 1 0 2 <z,z,z,z> 2 1 1 /
\endpicture} at 5.5 0 

\put{and} at 2.5 0 

\endpicture}
$$

A {\it valued translation quiver} $\Gamma = (\Gamma_0,\Gamma_1,\tau,v)$ is 
a translation quiver with two functions
$v',v''\:\Gamma_1 \to \Bbb N_1$ such that 
$$
 v'(\tau z,y)  = v''(y,z) \quad \text{and} \quad 
  v''(\tau z,y)  = v'(y,z),
$$
for any arrow $y \to z$ in $\Gamma$, where $z$ is a non-projective vertex;
we write $v = (v',v'')$ and call $v$ the {\it valuation} of $\Gamma$.
	\medskip
Dealing with meshes as exhibited above, it is sufficient to write down
the valuation for the arrows pointing south-east, say as follows:
$$
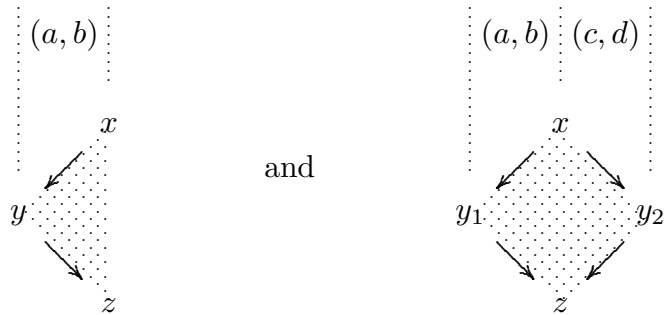
{\beginpicture
\setcoordinatesystem units <1cm,1cm>

\put{and} at 2.5 0 
\put{\beginpicture
\put{$x$} at 1 2
\put{$y$} at 0 1
\put{$z$} at 1 0
\arr{0.7 1.7}{0.3 1.3}
\arr{0.3  .7}{0.7  .3}
\put{$(a,b)$} at 0.5 3
\setdots <1mm>
\plot 0 3.3  0 1.5 /
\plot 1 3.3  1 2.5 /
\setshadegrid span <0.8mm>
\vshade 0 1 1 <z,z,z,z> 1 0 2 /
\endpicture} at 0 0 
\put{\beginpicture
\put{$x$} at 1 2
\put{$y_1$} at 0 1
\put{$y_2$} at 2 1
\put{$z$} at 1 0
\arr{0.7 1.7}{0.3 1.3}
\arr{1.3 1.7}{1.7 1.3}
\arr{0.3  .7}{0.7  .3}
\arr{1.7 .7}{1.3 .3}
\put{$(a,b)$} at 0.5 3
\put{$(c,d)$} at 1.5 3
\setdots <1mm>
\plot 0 3.3  0 1.5 /
\plot 1 3.3  1 2.5 /
\plot 2 3.3  2 1.5 /
\setshadegrid span <0.8mm>
\vshade 0 1 1 <z,z,z,z> 1 0 2  <z,z,z,z> 2 1 1 /
\endpicture} at 5.5 0 
\endpicture}
$$
this means on the right
that $v'(y_1,z) = a, v''(y_1,z) = b$ (thus $v'(x,y_1) = b, v''(x,y_1) = a$),
and that 
$v'(x,y_2) = c, v''(x,y_2) = d$ (thus $v'(y_2,z) = d, v''(y_2,z) = c$),
and similarly, on the left, that  
$v'(y,z) = a, v''(y,z) = b$ and $v'(x,y) = b, v''(x,y) = a$.
	\medskip
A function $g\:\Gamma_0 \to \Bbb Z$ is called 
{\it additive} provided
$$
 g(z) +g(\tau z) = \sum\nolimits_{y\in z^-} v'(y,z)g(y)
$$
for all non-projective vertices $z$ of $\Gamma$ (this is said to be the
corresponding {\it mesh relation}).
Thus, looking again at the meshes displayed above, we must have, on the left
$$
 g(z) = - g(x) + ag(y_1) + dg(y_2),
$$
and, on the right
$$
 g(z) = - g(x) + ag(y).
$$
	\medskip
The even-index triangle which we have presented above is 
the following translation quiver $\Gamma^{\ev}$:
its vertices are the pairs $(i,t)$
with $i,t\in \Bbb N_0$ such that $2i \le t$, with arrows $(i,t) \to (i,t+1)$
(we draw them as pointing south-west)
and, for $2i < t$, with arrows $(i,t) \to (i+1,t+1)$ (drawn as pointing south-east),
such that $\tau(i,t) = (i-1,t-2)$
provided $i\ge 1$, whereas the vertices of the form $(0,t)$ are projective. 
The vertices of the form $(i,2i)$ are said to lie on the
{\it pylon.}

We use the following valuation: the valuation of the arrows ending on the pylon
will be $(3,1),$ all other south-east arrows have valuation $(2,1).$ Here is
part of the triangle, the pylon being marked as a black line; the upper
row records the valuation:
$$
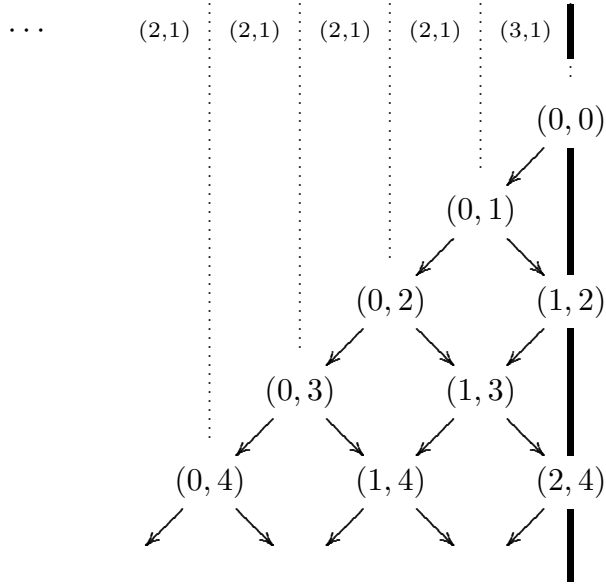
{\beginpicture
\setcoordinatesystem units <1cm,1cm>
\put{$(0,0)$} at 0 0 
\put{$(0,1)$} at -1 -1
\put{$(0,2)$} at -2 -2
\put{$(1,2)$} at 0 -2
\put{$(0,3)$} at -3 -3
\put{$(1,3)$} at -1 -3
\put{$(0,4)$} at -4 -4
\put{$(1,4)$} at -2 -4
\put{$(2,4)$} at 0 -4 
\arr{-.3 -.3}{-.7 -.7}

\arr{-1.3 -1.3}{-1.7 -1.7}
\arr{-.7 -1.3}{-.3 -1.7}

\arr{-2.3 -2.3}{-2.7 -2.7}
\arr{-1.7 -2.3}{-1.3 -2.7}
\arr{-.3 -2.3}{-.7 -2.7}

\arr{-3.3 -3.3}{-3.7 -3.7}
\arr{-2.7 -3.3}{-2.3 -3.7}
\arr{-1.3 -3.3}{-1.7 -3.7}
\arr{-.7 -3.3}{-.3 -3.7}

\arr{-4.3 -4.3}{-4.7 -4.7}
\arr{-3.7 -4.3}{-3.3 -4.7}
\arr{-2.3 -4.3}{-2.7 -4.7}
\arr{-1.7 -4.3}{-1.3 -4.7}
\arr{-.3 -4.3}{-.7 -4.7}

\put{$\ssize (3,1)$} at  -.5  1 
\multiput{$\ssize (2,1)$} at -1.5  1
 -2.5  1  -3.5  1  -4.5  1  /

\multiput{$\cdots$} at  -6 1 /

\setdots <1mm>
\plot 0 1.3  0 .5 /
\plot -1 1.3  -1 -.5 /
\plot -2 1.3  -2 -1.5 /
\plot -3 1.3  -3 -2.5 /
\plot -4 1.3  -4 -3.5 /

\put{} at 0 -5

\setsolid

\linethickness.7mm
\putrule from 0 -.3 to 0 -1.7
\putrule from 0 -2.3 to 0 -3.7
\putrule from 0 -4.3 to 0 -5.1

\putrule from 0 1.3 to 0 0.7

\endpicture}
$$
Note that $\Gamma^{\ev}$ is a subquiver of the 
valued translation quiver $\Bbb Z\Delta^{\ev}$, where $\Delta^{\ev}$ is the 
valued graph
$$
{\beginpicture
\setcoordinatesystem units <1.5cm,.7cm>
\multiput{$\circ$} at -3 0  -2 0  -1 0  0 0 /
\put{$-3$} at -3 0.5 
\put{$-2$} at -2 0.5 
\put{$-1$} at -1 0.5 
\put{$0$} at 0 0.5 
\plot -.8 0  -.2 0 /
\plot -.8 .15  -.2 .15 /
\plot -.8 -.15  -.2 -.15 /

\plot -2.8 .1  -2.2 .1 /
\plot -1.8 .1  -1.2 .1 /
 
\plot -2.8 -.1  -2.2 -.1 /
\plot -1.8 -.1  -1.2 -.1 /
\plot -2.6 .2 -2.4 0  -2.6 -.2 / 
\plot -1.6 .2 -1.4 0  -1.6 -.2 /

\plot -.65 .25 -.4 0  -.65 -.25 /

\multiput{$\dots$} at -4 0  /
\endpicture}
$$
(for the construction of $\Bbb Z\Delta$, where $\Delta$ is a
valued quiver, see [HPR]). If necessary, then we consider the underlying 
translation quiver
of $\Gamma^{\ev}$ and of $\Gamma^{\odd}$
(a valued translation quiver which will
be defined in section 2) as a subquiver of the translation quiver 
$\Bbb Z \Bbb A_\infty^\infty$; the latter quiver has as vertex set the set 
$\Bbb Z\times \Bbb Z$, there are arrows $(a,b) \to (a,b+1)$ and
$(a,b) \to (a+1,b+1)$, and the translation is given by $(a,b) \mapsto 
(a-1,b-2)$, for all $a,b\in \Bbb Z.$ Actually, for our considerations
it always will be sufficient to deal with the following subquiver of
$\Bbb Z \Bbb A_\infty^\infty$:

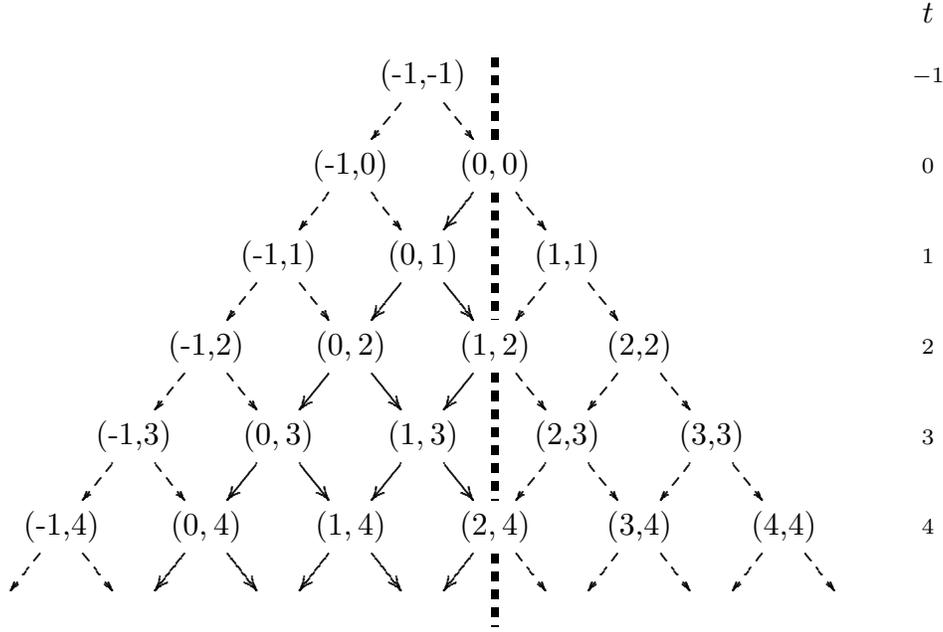
$$
{\beginpicture
\setcoordinatesystem units <.8cm,1cm>
\put{$(0,0)$} at 0 0 
\put{$(0,1)$} at -1 -1
\put{$(0,2)$} at -2 -2
\put{$(1,2)$} at 0 -2
\put{$(0,3)$} at -3 -3
\put{$(1,3)$} at -1 -3
\put{$(0,4)$} at -4 -4
\put{$(1,4)$} at -2 -4
\put{$(2,4)$} at 0 -4

\arr{-.3 -.3}{-.7 -.7}

\arr{-1.3 -1.3}{-1.7 -1.7}
\arr{-.7 -1.3}{-.3 -1.7}

\arr{-2.3 -2.3}{-2.7 -2.7}
\arr{-1.7 -2.3}{-1.3 -2.7}
\arr{-.3 -2.3}{-.7 -2.7}

\arr{-3.3 -3.3}{-3.7 -3.7}
\arr{-2.7 -3.3}{-2.3 -3.7}
\arr{-1.3 -3.3}{-1.7 -3.7}
\arr{-.7 -3.3}{-.3 -3.7}

\arr{-4.3 -4.3}{-4.7 -4.7}
\arr{-3.7 -4.3}{-3.3 -4.7}
\arr{-2.3 -4.3}{-2.7 -4.7}
\arr{-1.7 -4.3}{-1.3 -4.7}
\arr{-.3 -4.3}{-.7 -4.7}

\setdashes <1mm>
\arr{0.3 -0.3}{0.7 -0.7}
\arr{0.7 -1.3}{0.3 -1.7}
\arr{0.3 -2.3}{0.7 -2.7}
\arr{0.7 -3.3}{0.3 -3.7}
\arr{0.3 -4.3}{0.7 -4.7}

\arr{1.3 -1.3}{1.7 -1.7}
\arr{1.7 -2.3}{1.3 -2.7}
\arr{1.3 -3.3}{1.7 -3.7}
\arr{1.7 -4.3}{1.3 -4.7}

\arr{2.3 -2.3}{2.7 -2.7}
\arr{2.7 -3.3}{2.3 -3.7}
\arr{2.3 -4.3}{2.7 -4.7}

\arr{3.3 -3.3}{3.7 -3.7}
\arr{3.7 -4.3}{3.3 -4.7}

\arr{4.3 -4.3}{4.7 -4.7}

\arr{-1.3 0.7}{-1.7 0.3}
\arr{-0.7 0.7}{-0.3 0.3}
\arr{-2.3 -.3}{-2.7 -.7}
\arr{-1.7 -.3}{-1.3 -.7}
\arr{-3.3 -1.3}{-3.7 -1.7}
\arr{-2.7 -1.3}{-2.3 -1.7}
\arr{-4.3 -2.3}{-4.7 -2.7}
\arr{-3.7 -2.3}{-3.3 -2.7}
\arr{-5.3 -3.3}{-5.7 -3.7}
\arr{-4.7 -3.3}{-4.3 -3.7}
\arr{-6.3 -4.3}{-6.7 -4.7}
\arr{-5.7 -4.3}{-5.3 -4.7}

\put{(-1,-1)} at -1 1
\put{(-1,0)} at -2 0
\put{(-1,1)} at -3 -1
\put{(-1,2)} at -4 -2
\put{(-1,3)} at -5 -3
\put{(-1,4)} at -6 -4

\put{(1,1)} at 1 -1
\put{(2,2)} at 2 -2
\put{(2,3)} at 1 -3
\put{(3,3)} at 3 -3
\put{(3,4)} at 2 -4
\put{(4,4)} at 4 -4

\put{} at 0 -5

\setdashes <1mm>
\linethickness.7mm
\putrule from 0 1.2 to 0 .3
\putrule from 0 -.3 to 0 -1.7
\putrule from 0 -2.3 to 0 -3.7
\putrule from 0 -4.3 to 0 -5.1


\put{$\ssize -1$} at 6 1

\put{$\ssize 0$} at 6 0
\put{$\ssize 1$} at 6 -1
\put{$\ssize 2$} at 6 -2
\put{$\ssize 3$} at 6 -3
\put{$\ssize 4$} at 6 -4

\put{$t$} at 6 1.7

\endpicture}
$$
For better orientation, we have inserted the 
pylon for $\Gamma^{\ev}$ as a dotted vertical line (in 
the case of $\Gamma^{\odd}$ this will be the position of the 
left pylon).

	\medskip

If $g$ is an additive function on the translation quiver 
$\Gamma^{\ev}$, we write $g_i(t)$ instead
of $g(i,t)$. Note that any additive function $g$ on $\Gamma^{\ev}$ is
uniquely determined by the values $g(p)$ with $p$ projective,
thus by the values of $g_0.$ 
	\medskip
By definition, the function $d$ presented in the even-index triangle is an
{\it additive function on $\Gamma^{\ev}$.} And
{\it $d$ is the unique additive function on $\Gamma^{\ev}$ such that
$d_0(t) = 1$ for all $t \ge 0$} 
(thus with value $1$ on the projective vertices).
	\bigskip
{\bf Proposition 1.} {\it The function $g\:\Gamma^{\ev}_0 \to \Bbb Z$ 
is additive on $\Gamma^{\ev}$ if and only if it satisfies the following
hook condition for all $t\ge 1$:
$$
 g_i(t) = g_i(t-1) + \sum_{0\le j < i} g_j(t-i+j) 
$$
provided $2i < t$, and, for $t = 2i$:}
$$
 g_i(2i) = g_{i-1}(2i-1) + \sum_{0\le j < i} g_j(i+j).
$$
	\medskip
If we define $g_i(2i-1) = g_{i-1}(2i-1)$ (thus adding a column
on the right of the triangle, with vertex $(i,2i-1)$ in the row with index $2i-1$), 
then
the second condition in Proposition 1 has the same form as the 
first condition, with $t = 2i$, namely
$$
 g_i(2i) = g_i(2i-1) + \sum_{0\le j < i} g_j(2i-i+j),
$$ 
thus all conditions concern ``hooks'' in the even-index triangle 
as follows (the value at the circle being
obtained by adding the values at the bullets):

$$
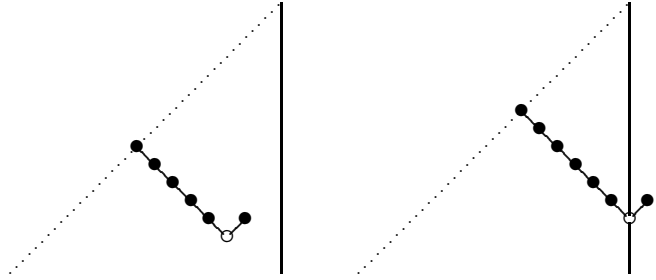
{\beginpicture
\setcoordinatesystem units <.2cm,.2cm>
\put{\beginpicture
\setdots <1mm>
\plot 0 0  15 15   /
\setsolid
\linethickness .2mm
\putrule from 15 0 to  15 15

\multiput{$\bullet$} at 7 7  8 6  9 5  10 4  11 3  13 3 /
\setsolid
\plot 7 7  11.8 2.2 / 
\plot 12.2 2.2 13 3 /

\put{$\circ$} at 12 2 

\multiput{} at 0 0  15 15 /

\endpicture} at 0 0
\put{\beginpicture
\setdots <1mm>
\plot -3 0  12 15  /
\multiput{$\bullet$} at 6 9  7 8  8 7  9 6  10 5  11 4  13 4 /
\setsolid
\plot 6 9  11.8 3.2 / 
\plot 12.2 3.2 13 4 /

\put{$\circ$} at 12 3

\multiput{} at -3 0  12 15 /
\setsolid
\linethickness .2mm
\putrule from 12 0 to  12 2.8
\putrule from 12 3.2 to  12 15

\endpicture} at 20 0
\endpicture}
$$

Here we should insert a remark concerning the difference between the
additivity property we encounter for the Fibonacci partition triangles and the
additivity property of the 
Pascal triangle. The additivity property for
the Pascal triangle means that any coefficient is the sum of the two upper
neighbors (one of them may be zero, if we are on the boundary). The hook condition means that we 
have to add not only the values of the two upper neighbors, but that we have to deal with
the values on a hook (but all the summands are still taken with multiplicity one). 
In contrast, the additivity property for a valued translation quiver (with arrows
pointing downwards) means that the value at the vertex $z$ is obtained by first taking
a certain 
linear combination of the values at the two upper neighbors (using positive coefficients
which may be different from $1$) and then 
subtracting the value at $\tau z$. 
Note that,
in general, such a mesh relation always involves subtraction, thus it
may lead to negative numbers. Of course, in our case, the equivalent hook condition
shows that we stay inside the set of positive integers.

	\medskip
Proof of Proposition 1, by induction on $t$. First, assume that $2i < t.$ Then the additivity
formula for $g$ yields:
$$
\align
 g_i(t) &= 2g_{i-1}(t-1) + g_i(t-1) - g_{i-1}(t-2) \cr
        &= g_i(t-1) + g_{i-1}(t-1) + (g_{i-1}(t-1)- g_{i-1}(t-2))  \cr
        &= g_i(t-1) + g_{i-1}(t-1) + \sum_{0\le j < i-1} g_j(t-i+j)  \cr
        &= g_i(t-1) + \sum_{0\le j < i} g_j(t-i+j).
\endalign
$$
Similarly, for $t = 2i > 0$ we get:
$$
\align
 g_i(2i) &= 3g_{i-1}(2i-1) - g_{i-1}(2i-2) \cr
        &= g_{i-1}(2i-1) + g_{i-1}(2i-1) + (g_{i-1}(2i-1)- g_{i-1}(2i-2))  \cr
        &= g_{i-1}(2i-1) + g_{i-1}(2i-1) + \sum_{0\le j < i-1} g_j(i+j)  \cr
        &= g_{i-1}(2i-1) + \sum_{0\le j < i} g_j(i+j).
\endalign
$$

The converse is shown in the same way.
	\bigskip
{\bf Corollary 1.} {\it The function $d_i(t)$ is a polynomial of degree $i$,
for any $i \ge 0,$} it is a monic linear combination of the binomial
coefficients $\binom t n$.
	\medskip
Proof: We use induction on $i$. For $i=0$, we deal with the constant polynomial
$d_0(t) = 1 = \binom t 0.$ Now let $i > 0.$ Then we have
$$
 d_i(t) - d_i(t-1) = \sum_{0\le j < i} d_j(t-i+j),
$$
and the right side is by induction a monic 
linear combination of the binomial
coefficients $\binom t 0, \binom t 1, \dots \binom t {i-1}$.	
Thus $g_i$ has to be a monic linear combination of the binomial
coefficients $\binom t 0, \binom t 1, \dots \binom t i$. 
In particular, $d_i(t)$ is a polynomial of degree $i$. 
	\bigskip

Here are the first polynomials $d_i(t)$:
$$
\align
 d_0(t) &= 1 = \binom t 0,\cr
 d_1(t) &= t = \binom t1, \cr
 d_2(t) &= \tfrac12 (t^2-t-6) = \binom t2- 3\binom t0, \cr
 d_3(t) &= \tfrac 16(t^3+3t^2-22t-18) = \binom t3+2\binom t2 -3\binom t1 - 3\binom t 0.
\endalign
$$

\noindent
It would be interesting to know a general formula for the polynomials $d_i(t)$.
	\bigskip

Looking back at the page showing the even-index triangle, we still have to
explain the entries $P_i$ in the left column. This concerns the 
Fibonacci modules $P_i = P_i(x)$  considered in [FR1] and [FR2], these are
indecomposable representations of the quiver obtained from the $3$-regular tree 
by choosing a bipartite orientation (if $i$ is even, $x$ has to be a sink, otherwise a source; in the following picture on the left, $x$ is
the vertex in the center). 
For example,
the top of the Fibonacci module $P_4$ has length $f_{8} = 21$ and its socle
has length $f_{10}= 55.$ The dimension vector of $P_4$ is of the form
as shown below on the left; the corresponding two rows $t=3$ and $t=4$ 
of the triangle (copied here on the right)
display one of the many walks in the support of $P_4$ which goes from the boundary to the center $x$.

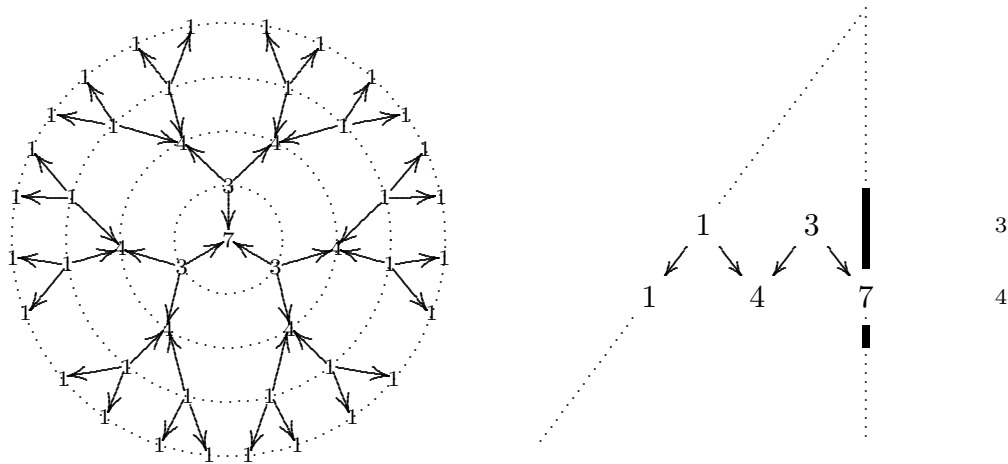
$$
{\beginpicture
\setcoordinatesystem units <0.6cm,0.6cm>
\multiput{\beginpicture

\put{$\ssize 7$} at 0 0
\multiput{$\ssize 3$} at -.866 -.5
  .866 -.5  0 1 /

\arrow <2mm> [0.25,0.75] from 0 0.9 to  0 0.2
\arrow <2mm> [0.25,0.75] from -.78 -.45
 to   -0.1 -0.05 
\arrow <2mm> [0.25,0.75] from .78 -.45
 to   0.1 -0.05 

\multiput{$\ssize 4$} at -.866 1.8 .866 1.8 /
\arrow <2mm> [0.25,0.75] from -0.1 1.1  to  -.78 1.75
\arrow <2mm> [0.25,0.75] from  0.1 1.1  to   .78 1.75
\multiput{$\ssize 1$} at -2.13 2.1  -1.1 2.8 
  2.13 2.1  1.1 2.8 /
\arrow <2mm> [0.25,0.75] from -2. 2.1     to  -.95 1.8 
\arrow <2mm> [0.25,0.75] from -1.07 2.7   to  -.866 1.9
\arrow <2mm> [0.25,0.75] from 2. 2.1      to  .95 1.8 
\arrow <2mm> [0.25,0.75] from 1.07 2.7    to  .866 1.9

\multiput{$\ssize 1$} at -3.27 2.3  -2.65 3 -1.7 3.62  -.7 3.93  /
\multiput{$\ssize 1$} at  3.27 2.3  2.65 3  1.7 3.62  .7 3.93 /

\arrow <2mm> [0.25,0.75] from -2.3 2.14   to  -3.15 2.3 
\arrow <2mm> [0.25,0.75] from -2.15 2.2    to  -2.6 2.9
\arrow <2mm> [0.25,0.75] from -1.15 2.9    to  -1.65 3.52
\arrow <2mm> [0.25,0.75] from -1.05 2.9    to  -.7 3.83 

\arrow <2mm> [0.25,0.75] from 2.3 2.14   to  3.15 2.3 
\arrow <2mm> [0.25,0.75] from 2.15 2.2    to  2.6 2.9
\arrow <2mm> [0.25,0.75] from 1.15 2.9    to  1.65 3.52
\arrow <2mm> [0.25,0.75] from 1.05 2.9    to  .7 3.83 

\startrotation by -0.5 0.866 about 0 0
\multiput{$\ssize 4$} at -.866 1.8 .866 1.8 /
\arrow <2mm> [0.25,0.75] from -0.1 1.1  to  -.78 1.75
\arrow <2mm> [0.25,0.75] from  0.1 1.1  to   .78 1.75
\multiput{$\ssize 1$} at -2.13 2.1  -1.1 2.8 
  2.13 2.1  1.1 2.8 /
\arrow <2mm> [0.25,0.75] from -2. 2.1     to  -.95 1.8 
\arrow <2mm> [0.25,0.75] from -1.07 2.7   to  -.866 1.9
\arrow <2mm> [0.25,0.75] from 2. 2.1      to  .95 1.8 
\arrow <2mm> [0.25,0.75] from 1.07 2.7    to  .866 1.9

\multiput{$\ssize 1$} at -3.27 2.3  -2.65 3 -1.7 3.62  -.7 3.93  /
\multiput{$\ssize 1$} at  3.27 2.3  2.65 3  1.7 3.62  .7 3.93 /

\arrow <2mm> [0.25,0.75] from -2.3 2.14   to  -3.15 2.3 
\arrow <2mm> [0.25,0.75] from -2.15 2.2    to  -2.6 2.9
\arrow <2mm> [0.25,0.75] from -1.15 2.9    to  -1.65 3.52
\arrow <2mm> [0.25,0.75] from -1.05 2.9    to  -.7 3.83 

\arrow <2mm> [0.25,0.75] from 2.3 2.14   to  3.15 2.3 
\arrow <2mm> [0.25,0.75] from 2.15 2.2    to  2.6 2.9
\arrow <2mm> [0.25,0.75] from 1.15 2.9    to  1.65 3.52
\arrow <2mm> [0.25,0.75] from 1.05 2.9    to  .7 3.83 
\stoprotation

\startrotation by -0.5 -0.866 about 0 0
\multiput{$\ssize 4$} at -.866 1.8 .866 1.8 /
\arrow <2mm> [0.25,0.75] from -0.1 1.1  to  -.78 1.75
\arrow <2mm> [0.25,0.75] from  0.1 1.1  to   .78 1.75
\multiput{$\ssize 1$} at -2.13 2.1  -1.1 2.8 
  2.13 2.1  1.1 2.8 /
\arrow <2mm> [0.25,0.75] from -2. 2.1     to  -.95 1.8 
\arrow <2mm> [0.25,0.75] from -1.07 2.7   to  -.866 1.9
\arrow <2mm> [0.25,0.75] from 2. 2.1      to  .95 1.8 
\arrow <2mm> [0.25,0.75] from 1.07 2.7    to  .866 1.9

\multiput{$\ssize 1$} at -3.27 2.3  -2.65 3 -1.7 3.62  -.7 3.93  /
\multiput{$\ssize 1$} at  3.27 2.3  2.65 3  1.7 3.62  .7 3.93 /

\arrow <2mm> [0.25,0.75] from -2.3 2.14   to  -3.15 2.3 
\arrow <2mm> [0.25,0.75] from -2.15 2.2    to  -2.6 2.9
\arrow <2mm> [0.25,0.75] from -1.15 2.9    to  -1.65 3.52
\arrow <2mm> [0.25,0.75] from -1.05 2.9    to  -.7 3.83 

\arrow <2mm> [0.25,0.75] from 2.3 2.14   to  3.15 2.3 
\arrow <2mm> [0.25,0.75] from 2.15 2.2    to  2.6 2.9
\arrow <2mm> [0.25,0.75] from 1.15 2.9    to  1.65 3.52
\arrow <2mm> [0.25,0.75] from 1.05 2.9    to  .7 3.83 
\stoprotation

\setdots <1mm>
\circulararc 360 degrees from 1 0 center at 0 0
\circulararc 360 degrees from 2 0 center at 0 0
\circulararc 360 degrees from 3 0 center at 0 0
\circulararc 360 degrees from 4 0 center at 0 0
\endpicture} at 0 0    /
\put{\beginpicture
\setcoordinatesystem units <.6cm,.8cm>
\put{$1$} at 0 0
\put{$1$} at 1 1
\put{$4$} at 2 0
\put{$3$} at 3 1
\put{$7$} at 4 0
\arr{0.7 0.7}{0.3 0.3}
\arr{1.3 0.7}{1.7 0.3}
\arr{2.7 0.7}{2.3 0.3}
\arr{3.3 0.7}{3.7 0.3}

\linethickness.7mm
\putrule from 4 -.7 to 4 -.4
\putrule from 4  .4 to 4 1.5

\put{$\ssize 3$} at 6.5 1
\put{$\ssize 4$} at 6.5 0

\setdots <1mm>
\plot 1.3 1.3  4 4  4 1.3 /
\plot -.3 -.3  -2 -2 /
\plot 4 -.3  4 -2 /
\endpicture} at 11 -.5

\endpicture}
$$
	\medskip
To be precise,
the numbers $d_i(t)$ are categorified by the modules $P_n = P_n(x)$ as follows:
they provide the Jordan-H\"older multiplicities of these modules.
First, let us consider the socle of $P_t$; the composition factors in
the socle of $P_t$ are of the form $S(z)$ with $0 \le D(x,z) \le t$ and
$D(x,z) \equiv t \mod 2$,
where, $D(x,z)$ denotes the distance between $x,z$ in the $3$-regular tree
(and $S(z)$ is the one-dimensional representation with $\dim S(z)_z \neq 0$).
There is the following multiplicity formula:
	\medskip
\Rahmen
{$
  d_i(t) = \dim P_t(x)_z \qquad \text{for}\quad D(x,z) = |t-2i|.
$}
	\medskip
(Remark: In case we assume, as we usually do, that $i \le \frac t2$, then we
may just write $D(x,z) = t-2i.$ However, it is sometimes convenient to
consider also the values $d_i(t)$ for $i > \frac t2$, where, by definition,
$d_i(t) = d_{t-i}(t)$. There is the second convention that $d_i(t) = 0$
for $i < 0$ and all $t$; again, the rule above is valid also for $i < 0$,
since $t-2i > t$ for $i < 0$, and $\dim P_t(x)_z = 0$
in case $D(x,z) > t.$)
	\medskip
Note that the socle of $P_t(x)$ is isomorphic to the top of $P_{t+1}(x)$, thus
$d_i(t)$ is also equal to the dimension of $\dim P_{t+1}(x)_z$ 
where $D(x,z) = |t-2i|.$
	\medskip
We should add the following warning: the module $P_n(x)$ with $n$ even
is only defined in case $x$ is a sink, for $n$ odd, if $x$ is a source;
in particular, the modules $P_t(x)$ and $P_{t+1}(x)$ are {\bf not} defined for 
the same quiver.

\vfill\eject

{\bf 2. The odd-index Fibonacci partition triangle}

$$
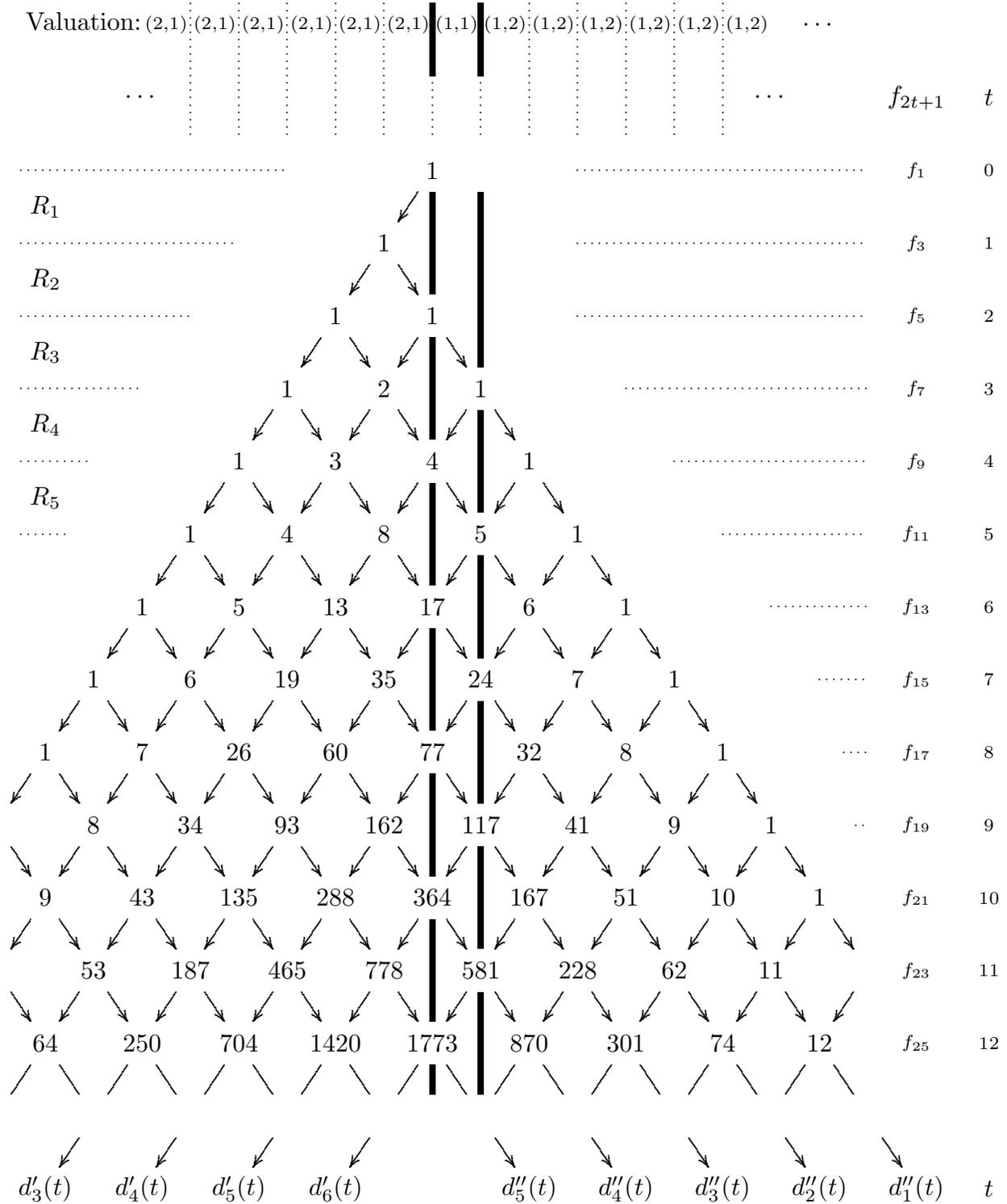
{\beginpicture
\setcoordinatesystem units <.665cm,1cm>

\multiput{$\ssize (2,1)$} at -1.5  1  
 -2.5  1  -3.5  1  -4.5  1  -5.5  1  -6.5  1   /
\put{$\ssize (1,1)$} at  -.5  1 

\multiput{$\ssize (1,2)$} at .5 1  1.5 1  2.5 1  3.5 1  4.5 1  5.5 1 /

\multiput{$\cdots$} at   7 1  -7 0  6 0 /

\put{Valuation:} [l] at -9.4 1.05

\setdots <1mm>
\plot 0 -.5  0 1.3 /
\plot -1 -.5  -1 1.3 /
\plot -2 -.5  -2 1.3 /
\plot -3 -.5  -3 1.3 /
\plot -4 -.5  -4 1.3 /
\plot -5 -.5  -5 1.3 /
\plot -6 -.5  -6 1.3 /

\plot 1 -.5  1 1.3 /
\plot 2 -.5  2 1.3 /
\plot 3 -.5  3 1.3 /
\plot 4 -.5  4 1.3 /
\plot 5 -.5  5 1.3 /

\setsolid

\multiput{$1$} at -1 -1 -2 -2  -3 -3  -4 -4  -5 -5  -6 -6  -7 -7  -8 -8  -9 -9 /

\put{$1$} at -1 -3 
\put{$2$} at -2 -4 
\put{$3$} at -3 -5 
\put{$4$} at -4 -6 
\put{$5$} at -5 -7 
\put{$6$} at -6 -8 
\put{$7$} at -7 -9 
\put{$8$} at -8 -10 
\put{$9$} at -9 -11 

\put{$1$} at 0 -4
\put{$4$} at -1 -5
\put{$8$} at -2 -6 
\put{$13$} at -3 -7 
\put{$19$} at -4 -8 
\put{$26$} at -5 -9 
\put{$34$} at -6 -10 
\put{$43$} at -7 -11
\put{$53$} at -8 -12 

\put{$5$} at 0 -6
\put{$17$} at -1 -7 
\put{$35$} at -2 -8 
\put{$60$} at -3 -9 
\put{$93$} at -4 -10 

\put{$77$} at -1 -9 
\put{$162$} at -2 -10

\arr{-1.3 -1.3}{-1.7 -1.7}

\arr{-2.3 -2.3}{-2.7 -2.7}
\arr{-1.7 -2.3}{-1.3 -2.7}

\arr{-3.3 -3.3}{-3.7 -3.7}
\arr{-2.7 -3.3}{-2.3 -3.7}
\arr{-1.3 -3.3}{-1.7 -3.7}
\arr{-.7 -3.3}{-.3 -3.7}

\arr{-4.3 -4.3}{-4.7 -4.7}
\arr{-3.7 -4.3}{-3.3 -4.7}
\arr{-2.3 -4.3}{-2.7 -4.7}
\arr{-1.7 -4.3}{-1.3 -4.7}
\arr{-.3 -4.3}{-.7 -4.7}

\arr{-5.3 -5.3}{-5.7 -5.7}
\arr{-4.7 -5.3}{-4.3 -5.7}
\arr{-3.3 -5.3}{-3.7 -5.7}
\arr{-2.7 -5.3}{-2.3 -5.7}
\arr{-1.3 -5.3}{-1.7 -5.7}
\arr{-.7 -5.3}{-.3 -5.7}

\arr{-6.3 -6.3}{-6.7 -6.7}
\arr{-5.7 -6.3}{-5.3 -6.7}
\arr{-4.3 -6.3}{-4.7 -6.7}
\arr{-3.7 -6.3}{-3.3 -6.7}
\arr{-2.3 -6.3}{-2.7 -6.7}
\arr{-1.7 -6.3}{-1.3 -6.7}
\arr{-.3 -6.3}{-.7 -6.7}

\arr{-7.3 -7.3}{-7.7 -7.7}
\arr{-6.7 -7.3}{-6.3 -7.7}
\arr{-5.3 -7.3}{-5.7 -7.7}
\arr{-4.7 -7.3}{-4.3 -7.7}
\arr{-3.3 -7.3}{-3.7 -7.7}
\arr{-2.7 -7.3}{-2.3 -7.7}
\arr{-1.3 -7.3}{-1.7 -7.7}
\arr{-.7 -7.3}{-.3 -7.7}

\arr{-8.3 -8.3}{-8.7 -8.7}
\arr{-7.7 -8.3}{-7.3 -8.7}
\arr{-6.3 -8.3}{-6.7 -8.7}
\arr{-5.7 -8.3}{-5.3 -8.7}
\arr{-4.3 -8.3}{-4.7 -8.7}
\arr{-3.7 -8.3}{-3.3 -8.7}
\arr{-2.3 -8.3}{-2.7 -8.7}
\arr{-1.7 -8.3}{-1.3 -8.7}
\arr{-.3 -8.3}{-.7 -8.7}

\arr{-9.3 -9.3}{-9.7 -9.7}
\arr{-8.7 -9.3}{-8.3 -9.7}
\arr{-7.3 -9.3}{-7.7 -9.7}
\arr{-6.7 -9.3}{-6.3 -9.7}
\arr{-5.3 -9.3}{-5.7 -9.7}
\arr{-4.7 -9.3}{-4.3 -9.7}
\arr{-3.3 -9.3}{-3.7 -9.7}
\arr{-2.7 -9.3}{-2.3 -9.7}
\arr{-1.3 -9.3}{-1.7 -9.7}
\arr{-.7 -9.3}{-.3 -9.7}

\arr{-9.7 -10.3}{-9.3 -10.7}
\arr{-8.3 -10.3}{-8.7 -10.7}
\arr{-7.7 -10.3}{-7.3 -10.7}
\arr{-6.3 -10.3}{-6.7 -10.7}
\arr{-5.7 -10.3}{-5.3 -10.7}
\arr{-4.3 -10.3}{-4.7 -10.7}
\arr{-3.7 -10.3}{-3.3 -10.7}
\arr{-2.3 -10.3}{-2.7 -10.7}
\arr{-1.7 -10.3}{-1.3 -10.7}
\arr{-.3 -10.3}{-.7 -10.7}

\arr{-9.3 -11.3}{-9.7 -11.7}
\arr{-8.7 -11.3}{-8.3 -11.7}
\arr{-7.3 -11.3}{-7.7 -11.7}
\arr{-6.7 -11.3}{-6.3 -11.7}
\arr{-5.3 -11.3}{-5.7 -11.7}
\arr{-4.7 -11.3}{-4.3 -11.7}
\arr{-3.3 -11.3}{-3.7 -11.7}
\arr{-2.7 -11.3}{-2.3 -11.7}
\arr{-1.3 -11.3}{-1.7 -11.7}
\arr{-.7 -11.3}{-.3 -11.7}

\arr{-9.7 -12.3}{-9.3 -12.7}
\arr{-8.3 -12.3}{-8.7 -12.7}
\arr{-7.7 -12.3}{-7.3 -12.7}
\arr{-6.3 -12.3}{-6.7 -12.7}
\arr{-5.7 -12.3}{-5.3 -12.7}
\arr{-4.3 -12.3}{-4.7 -12.7}
\arr{-3.7 -12.3}{-3.3 -12.7}
\arr{-2.3 -12.3}{-2.7 -12.7}
\arr{-1.7 -12.3}{-1.3 -12.7}
\arr{-.3 -12.3}{-.7 -12.7}

\setsolid

\plot -9.3 -13.3 -9.7 -13.7 /
\plot -8.7 -13.3 -8.3 -13.7 /
\plot -7.3 -13.3 -7.7 -13.7 /
\plot -6.7 -13.3 -6.3 -13.7 /
\plot -5.3 -13.3 -5.7 -13.7 /
\plot -4.7 -13.3 -4.3 -13.7 /
\plot -3.3 -13.3 -3.7 -13.7 /
\plot -2.7 -13.3 -2.3 -13.7 /
\plot -1.3 -13.3 -1.7 -13.7 /
\plot  -.7 -13.3  -.3 -13.7 /

\linethickness.7mm
\putrule from 0 -1.3 to 0 -3.7
\putrule from 0 -4.3 to 0 -5.7
\putrule from 0 -6.3 to 0 -7.7
\putrule from 0 -8.3 to 0 -9.7
\putrule from 0 -10.3 to 0 -11.7
\putrule from 0 -12.3 to 0 -13.7

\putrule from -1 -1.3 to -1 -2.7
\putrule from -1 -3.3 to -1 -4.7
\putrule from -1 -5.3 to -1 -6.7
\putrule from -1 -7.3 to -1 -8.7
\putrule from -1 -9.3 to -1 -10.7
\putrule from -1 -11.3 to -1 -12.7
\putrule from -1 -13.3 to -1 -13.7

\setsolid

\putrule from 0 1.3 to 0 .3
\putrule from -1 1.3 to -1 .3

\put{$\ssize 0$} at 10.5 -1
\put{$\ssize 1$} at 10.5 -2
\put{$\ssize 2$} at 10.5 -3
\put{$\ssize 3$} at 10.5 -4
\put{$\ssize 4$} at 10.5 -5
\put{$\ssize 5$} at 10.5 -6
\put{$\ssize 6$} at 10.5 -7
\put{$\ssize 7$} at 10.5 -8
\put{$\ssize 8$} at 10.5 -9
\put{$\ssize 9$} at 10.5 -10
\put{$\ssize 10$} at 10.5 -11
\put{$\ssize 11$} at 10.5 -12
\put{$\ssize 12$} at 10.5 -13

\put{$d'_3(t)$} at -9 -15
\put{$d'_4(t)$} at -7 -15
\put{$d'_5(t)$} at -5 -15
\put{$d'_6(t)$} at -3 -15

\put{$t$} at 10.5 -15

\put{$d''_1(t)$} at 9 -15
\put{$d''_2(t)$} at 7 -15
\put{$d''_3(t)$} at 5 -15
\put{$d''_4(t)$} at 3 -15
\put{$d''_5(t)$} at 1 -15

\arr{-8.3 -14.3}{-8.7 -14.7}
\arr{-6.3 -14.3}{-6.7 -14.7}
\arr{-4.3 -14.3}{-4.7 -14.7}
\arr{-2.3 -14.3}{-2.7 -14.7}

\arr{8.3 -14.3}{8.7 -14.7}
\arr{6.3 -14.3}{6.7 -14.7}
\arr{4.3 -14.3}{4.7 -14.7}
\arr{2.3 -14.3}{2.7 -14.7}
\arr{0.3 -14.3}{.7 -14.7}

\arr{.3 -4.3}{.7 -4.7}

\arr{1.3 -5.3}{1.7 -5.7}
\arr{.7 -5.3}{.3 -5.7}

\arr{2.3 -6.3}{2.7 -6.7}
\arr{1.7 -6.3}{1.3 -6.7}
\arr{.3 -6.3}{.7 -6.7}

\arr{3.3 -7.3}{3.7 -7.7}
\arr{2.7 -7.3}{2.3 -7.7}
\arr{1.3 -7.3}{1.7 -7.7}
\arr{.7 -7.3}{.3 -7.7}

\arr{4.3 -8.3}{4.7 -8.7}
\arr{3.7 -8.3}{3.3 -8.7}
\arr{2.3 -8.3}{2.7 -8.7}
\arr{1.7 -8.3}{1.3 -8.7}
\arr{.3 -8.3}{.7 -8.7}

\arr{5.3 -9.3}{5.7 -9.7}
\arr{4.7 -9.3}{4.3 -9.7}
\arr{3.3 -9.3}{3.7 -9.7}
\arr{2.7 -9.3}{2.3 -9.7}
\arr{1.3 -9.3}{1.7 -9.7}
\arr{.7 -9.3}{.3 -9.7}

\arr{6.3 -10.3}{6.7 -10.7}
\arr{5.7 -10.3}{5.3 -10.7}
\arr{4.3 -10.3}{4.7 -10.7}
\arr{3.7 -10.3}{3.3 -10.7}
\arr{2.3 -10.3}{2.7 -10.7}
\arr{1.7 -10.3}{1.3 -10.7}
\arr{.3 -10.3}{.7 -10.7}

\arr{7.3 -11.3}{7.7 -11.7}
\arr{6.7 -11.3}{6.3 -11.7}
\arr{5.3 -11.3}{5.7 -11.7}
\arr{4.7 -11.3}{4.3 -11.7}
\arr{3.3 -11.3}{3.7 -11.7}
\arr{2.7 -11.3}{2.3 -11.7}
\arr{1.3 -11.3}{1.7 -11.7}
\arr{.7 -11.3}{.3 -11.7}

\arr{7.7 -12.3}{7.3 -12.7}
\arr{6.3 -12.3}{6.7 -12.7}
\arr{5.7 -12.3}{5.3 -12.7}
\arr{4.3 -12.3}{4.7 -12.7}
\arr{3.7 -12.3}{3.3 -12.7}
\arr{2.3 -12.3}{2.7 -12.7}
\arr{1.7 -12.3}{1.3 -12.7}
\arr{.3 -12.3}{.7 -12.7}

\plot 7.3 -13.3 7.7 -13.7 /
\plot 6.7 -13.3 6.3 -13.7 /
\plot 5.3 -13.3 5.7 -13.7 /
\plot 4.7 -13.3 4.3 -13.7 /
\plot 3.3 -13.3 3.7 -13.7 /
\plot 2.7 -13.3 2.3 -13.7 /
\plot 1.3 -13.3 1.7 -13.7 /
\plot .7 -13.3  .3 -13.7 /

\put{$1$} at 1 -5  
\put{$1$} at 2 -6  
\put{$1$} at 3 -7  
\put{$1$} at 4 -8  
\put{$1$} at 5 -9  
\put{$1$} at 6 -10  
\put{$1$} at 7 -11

\put{$6$} at 1 -7  
\put{$7$} at 2 -8
\put{$8$} at 3 -9  
\put{$9$} at 4 -10  
\put{$10$} at 5 -11  
\put{$11$} at 6 -12

\put{$24$} at 0 -8  
\put{$32$} at 1 -9  
\put{$117$} at 0 -10
\put{$41$} at 2 -10  
\put{$51$} at 3 -11  
\put{$62$} at 4 -12  
\put{$135$} at -5 -11
\put{$288$} at -3 -11
\put{$364$} at -1 -11
\put{$167$} at 1 -11
\put{$187$} at -6 -12
\put{$465$} at -4 -12
\put{$778$} at -2 -12
\put{$581$} at 0 -12
\put{$228$} at 2 -12

\put{$64$} at -9 -13
\put{$250$} at -7 -13
\put{$704$} at -5 -13
\put{$1420$} at -3 -13
\put{$1773$} at -1 -13
\put{$870$} at 1 -13
\put{$301$} at 3 -13
\put{$74$} at 5 -13
\put{$12$} at 7 -13

\setdots <1mm>
\plot -9.5 -1  -4 -1 /
\plot -9.5 -2  -5 -2 /
\plot -9.5 -3  -6 -3 /
\plot -9.5 -4  -7 -4 /
\plot -9.5 -5  -8 -5 /
\plot -9.5 -6  -8.5 -6 /

\plot 2 -1  8 -1 /
\plot 2 -2  8 -2 /
\plot 2 -3  8 -3 /
\plot 3 -4  8 -4 /
\plot 4 -5  8 -5 /
\plot 5 -6  8 -6 /
\plot 6 -7  8 -7 /
\plot 7 -8  8 -8 /
\plot 7.5 -9  8 -9 /
\plot 7.75 -10  8 -10 /

\put{$R_1$} at -9 -1.5
\put{$R_2$} at -9 -2.5
\put{$R_3$} at -9 -3.5
\put{$R_4$} at -9 -4.5
\put{$R_5$} at -9 -5.5

\put{$f_{2t+1}$} at 9 0
\put{$t$} at 10.5 0

\put{$\ssize f_1$} at 9 -1
\put{$\ssize f_3$} at 9 -2
\put{$\ssize f_5$} at 9 -3
\put{$\ssize f_{7}$} at 9 -4
\put{$\ssize f_{9}$} at 9 -5
\put{$\ssize f_{11}$} at 9 -6
\put{$\ssize f_{13}$} at 9 -7
\put{$\ssize f_{15}$} at 9 -8
\put{$\ssize f_{17}$} at 9 -9
\put{$\ssize f_{19}$} at 9 -10
\put{$\ssize f_{21}$} at 9 -11
\put{$\ssize f_{23}$} at 9 -12
\put{$\ssize f_{25}$} at 9 -13

\endpicture}
$$

\vfill\eject

\noindent
Again, let us first explain the numbers displayed here. As before, the rows are 
indexed by $t = 0,1,2,\dots$ as shown on the right. Note that we deal here with
a triangle
only after removing the row $t=0$; alternatively, we may interprete the display
as being formed by two triangles, separated by the left pylon. 

The entries in the rows with index $0$ and $1$ are labeled $d'_0(0) = d''_{-1}(0)$ and $d'_0(1) = d''_0(1)$,
respectively. Consider now a row with index $t\ge 2$.
The entries will be labeled, from left to right, by
$d'_0(t), d'_1(t), \dots, d'_i(t), \dots, d'_{t-1}(t)$, 
and we put $d''_i(t) = d'_{t-i-1}(t)$ for $0 \le i \le t-1$. 
Usually, we will need the notation $d'$ for the numbers up to the second
pylon from the left, 
the notation $d''$ for those starting with the first pylon and going to
the right. For example,
for $t = 5,$ the entries from left to right are labeled 
$$
 d'_0(5) = 1,\quad d'_1(5) = 4, \quad d'_2(5)= 8,\quad
 d'_3(5) = d''_1(5) = 5, \quad d''_0(5) = 1.
$$
If necessary, we will write $d'_i(t) = 0$ for $i < 0$. Also, let
$d''_i(t) = 0$ for $i < 0$, unless $i = {-1}$ and $t = 0.$
Using the rule $d''_i(t) = d'_{t-i-1}(t)$ we have defined in this way $d'_i(t)$
and $d''_i(t)$ for all $i \in \Bbb Z.$ 

The entries $d'_i(t)$ and $d''_i(t)$ displayed 
are calculated from the numbers $u_i[j]$ of the paper [FR2] according to the rule
$$
\align
 d'_i(t) &= u_{\lceil t/2 \rceil}[t-2i],\cr
 d''_i(t) &= u_{\lceil t/2 \rceil}[-t+2+2i].
\endalign
$$

Using this notation, the partition formulae of [FR2] again can be written in
a unified form, as follows:
	\medskip
{\bf Partition formula for the Fibonacci numbers $f_n$ with odd index $n$:}
	\medskip
\Rahmen{$
 \sum\limits_{0\le i < t/2} 2^{t-2i}d'_i(t)+
 \sum\limits_{0\le i < (t-3)/2} 2^{t-2i}d''_i(t) \ =  f_{2t+1}
$}

\noindent
Thus, {\it the row $t$ of the triangle yields the Fibonacci number $f_{2t+1}$.}
As examples, let us look at the rows $t= 5$ and $t=6$, they yield
$$
\align
 (2^5\cdot 1 + 2^3\cdot 4 + 2^1\cdot 8) +
 (2^0\cdot 5 + 2^2\cdot 1) &=\;\; 89 = f_{11}    \tag{$t=5$}\cr
 (2^6\cdot 1 + 2^4\cdot 5 + 2^2\cdot 13 + 2^0\cdot 17) + 
 (2^1\cdot 6 + 2^3\cdot 1)  &= 233 = f_{13}  \tag{$t=6$}
\endalign
$$

	\bigskip
Let us describe in detail the translation quiver $\Gamma^{\odd}$ used here.
Its vertices are the pair $(0,0)$ as well as all the pairs $(i,t)$ 
with integers $0 \le i < t$; there are the south-west arrows $(i,t)$ to $(i,t+1)$
as well as the south-east arrows $(i,t) \to (i+1,t+1)$, where $(i,t) \neq (0,0)$.
All the vertices of the form $(0,t)$ and $(i,i+1)$ with $i\ge 2$ are
projective, and $\tau(i,t) = (i-1,t-2)$ for the remaining vertices.
Now we mark two pylons: the first one consists of the vertices of the form
$(i,2i)$, the second one of the vertices $(i,2i-1)$.
All the arrows with both starting vertex and terminal vertex on a pylon have
valuation $(1,1)$; the (south-east) arrows $(i,t) \to (i+1,t+1)$ with $2i < t$
have valuation $(2,1)$, the remaining south-east arrows have valuation $(1,2)$.
$$
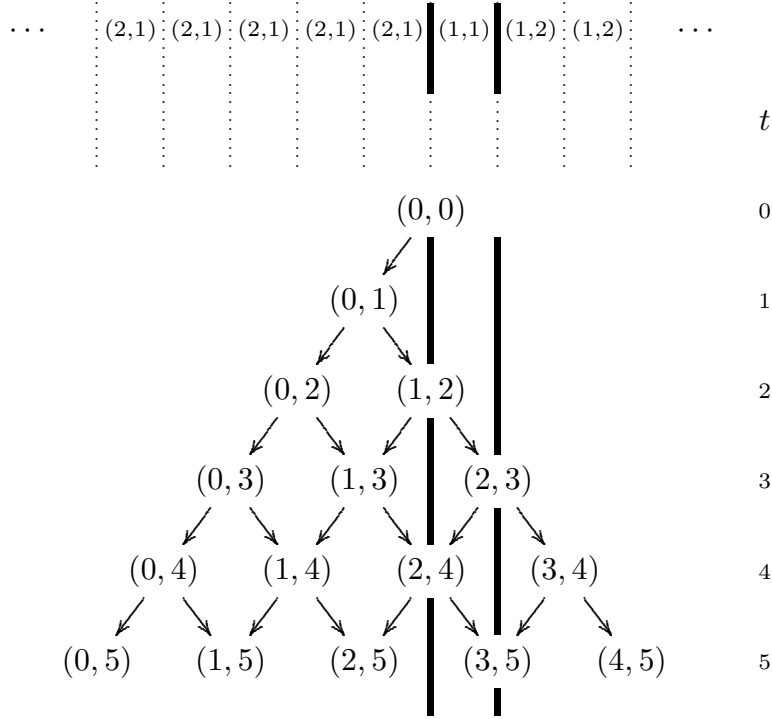
{\beginpicture
\setcoordinatesystem units <.74cm,1cm>

\multiput{$\ssize (2,1)$} at -1.5  1  
 -2.5  1  -3.5  1  -4.5  1  -5.5  1     /
\put{$\ssize (1,1)$} at  -.5  1 

\multiput{$\ssize (1,2)$} at .5 1  1.5 1   /

\multiput{$\cdots$} at   3 1   -7 1 /


\setdots <1mm>
\plot 0 -.5  0 1.3 /
\plot -1 -.5  -1 1.3 /
\plot -2 -.5  -2 1.3 /
\plot -3 -.5  -3 1.3 /
\plot -4 -.5  -4 1.3 /
\plot -5 -.5  -5 1.3 /
\plot -6 -.5  -6 1.3 /

\plot 1 -.5  1 1.3 /
\plot 2 -.5  2 1.3 /

\setsolid

\put{$(0,0)$} at -1 -1 
\put{$(0,1)$} at -2 -2 
\put{$(0,2)$} at -3 -3 
\put{$(0,3)$} at -4 -4 
\put{$(0,4)$} at -5 -5 
\put{$(0,5)$} at -6 -6

\put{$(1,2)$} at -1 -3 
\put{$(1,3)$} at -2 -4 
\put{$(1,4)$} at -3 -5 
\put{$(1,5)$} at -4 -6 

\put{$(2,3)$} at 0 -4
\put{$(2,4)$} at -1 -5
\put{$(2,5)$} at -2 -6 

\put{$(3,5)$} at 0 -6

\arr{-1.3 -1.3}{-1.7 -1.7}

\arr{-2.3 -2.3}{-2.7 -2.7}
\arr{-1.7 -2.3}{-1.3 -2.7}

\arr{-3.3 -3.3}{-3.7 -3.7}
\arr{-2.7 -3.3}{-2.3 -3.7}
\arr{-1.3 -3.3}{-1.7 -3.7}
\arr{-.7 -3.3}{-.3 -3.7}

\arr{-4.3 -4.3}{-4.7 -4.7}
\arr{-3.7 -4.3}{-3.3 -4.7}
\arr{-2.3 -4.3}{-2.7 -4.7}
\arr{-1.7 -4.3}{-1.3 -4.7}
\arr{-.3 -4.3}{-.7 -4.7}

\arr{-5.3 -5.3}{-5.7 -5.7}
\arr{-4.7 -5.3}{-4.3 -5.7}
\arr{-3.3 -5.3}{-3.7 -5.7}
\arr{-2.7 -5.3}{-2.3 -5.7}
\arr{-1.3 -5.3}{-1.7 -5.7}
\arr{-.7 -5.3}{-.3 -5.7}

\setsolid

\linethickness.7mm
\putrule from 0 -1.3 to 0 -3.7
\putrule from 0 -4.3 to 0 -5.7
\putrule from 0 -6.3 to 0 -6.6

\putrule from -1 -1.3 to -1 -2.7
\putrule from -1 -3.3 to -1 -4.7
\putrule from -1 -5.3 to -1 -6.6

\putrule from 0 1.3 to 0 .3
\putrule from -1 1.3 to -1 .3

\put{$\ssize 0$} at 4 -1
\put{$\ssize 1$} at 4 -2
\put{$\ssize 2$} at 4 -3
\put{$\ssize 3$} at 4 -4
\put{$\ssize 4$} at 4 -5
\put{$\ssize 5$} at 4 -6

\put{$t$} at 4 0

\arr{.3 -4.3}{.7 -4.7}

\arr{1.3 -5.3}{1.7 -5.7}
\arr{.7 -5.3}{.3 -5.7}

\put{$(3,4)$} at 1 -5  
\put{$(4,5)$} at 2 -6  

\endpicture}
$$
(Let us stress that this pair of pylons indicates completely different
mesh rules than the single pylon used in the even-index triangle.)

	\medskip
{\it The function $d'$ is an additive function on $\Gamma^{\odd}$,
this is the unique additive function on $\Gamma^{\odd}$ such that
$d'_0(t) = 1$ for all $t \ge 0$ and $d''_0(t) = 1$ for all $t \ge 3$}
(thus with values $1$ on the projective vertices).
	\medskip

Remark. Note that $\Gamma^{\odd}$ can be considered as a subquiver of the
valued translation quiver $\Bbb Z\Delta^{\odd}$, where $\Delta^{\odd}$ is the 
valued graph
$$
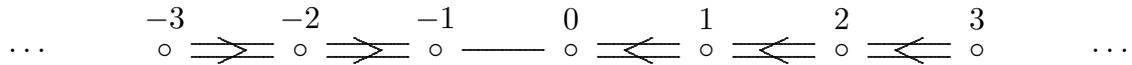
{\beginpicture
\setcoordinatesystem units <1.5cm,.7cm>
\multiput{$\circ$} at -3 0  -2 0  -1 0  0 0  1 0  2 0  3 0 /
\put{$-3$} at -3 0.5 
\put{$-2$} at -2 0.5 
\put{$-1$} at -1 0.5 
\put{$0$} at 0 0.5 
\put{$1$} at 1 0.5 
\put{$2$} at 2 0.5 
\put{$3$} at 3 0.5 
\plot -.8 0  -.2 0 /
\plot -2.8 .1  -2.2 .1 /
\plot -1.8 .1  -1.2 .1 /
\plot .8 .1  .2 .1 /
\plot 1.8 .1  1.2 .1 /
\plot 2.8 .1  2.2 .1 /
 
\plot -2.8 -.1  -2.2 -.1 /
\plot -1.8 -.1  -1.2 -.1 /
\plot .8 -.1  .2 -.1 /
\plot 1.8 -.1  1.2 -.1 /
\plot 2.8 -.1  2.2 -.1 /
\plot -2.6 .2 -2.4 0  -2.6 -.2 / 
\plot -1.6 .2 -1.4 0  -1.6 -.2 /
 \plot .6 .2  .4 0    .6 -.2 / 
\plot 1.6 .2  1.4 0  1.6 -.2 / 
\plot 2.6 .2  2.4 0  2.6 -.2 /

\multiput{$\dots$} at -4 0  4 0 /
\endpicture}
$$
as mentioned at the end of section 3 of [FR2].

	\medskip
For additive functions $g$ on $\Gamma^{\odd}$, 
we write $g'_i(t) = g(i,t)$, usually for $2i \le t+1$, 
and $g''_i(t) = g(t-1-i,t)$, usually for $2i \ge t.$
There is again a hook characterization of additivity:
	\bigskip
{\bf Proposition 2.} {\it The function $g\:\Gamma^{\odd}_0 \to \Bbb Z$ 
is additive on $\Gamma^{\odd}$ if and only if it satisfies the following
hook conditions:

\item{\rm(a)} If $2i \le t$, then 
$$
 g'_i(t) = g'_i(t-1) + \sum_{0\le j < i} g'_j(t-i+j).
$$

\item{\rm (b)} 
If $2i > t$, then} 
$$
 g''_i(t) = g''_i(t-1) + \sum_{0\le j < i} g''_j(t-i+j).
$$
	\medskip

The following pictures indicate the position of the hooks  
in the odd-index triangle (again, the value at the circle is
obtained by adding the values at the bullets):
$$
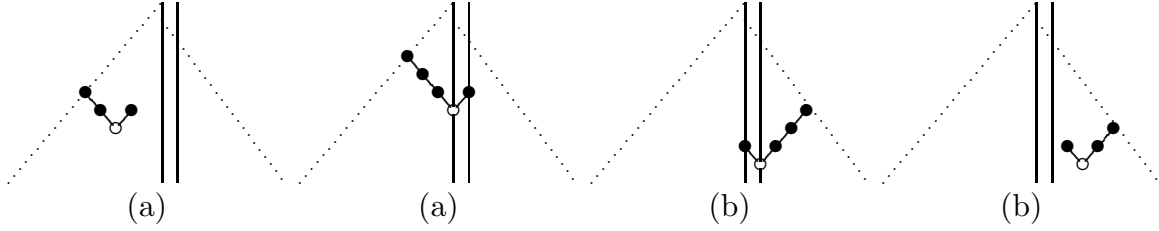
{\beginpicture
\setcoordinatesystem units <.17cm,.2cm>
\put{\beginpicture
\setdots <1mm>
\plot -1 -1  9 9 /
\plot 9 8  17 -1 /

\setsolid
\linethickness .2mm
\putrule from 9 -1 to  9 9
\putrule from 10 -1 to  10 9

\multiput{$\bullet$} at 4 4  5 3  7 3   /
\setsolid
\plot 4 4  5.8 2.2  / 
\plot 6.2 2.2 7 3 /

\put{$\circ$} at 6 2 

\multiput{} at 0 0  16 9 /
\put{(a)} at 8 -2.4
\endpicture} at 0 0
\put{\beginpicture
\setdots <1mm>
\plot -1 -1  9 9 /
\plot 9 8  17 -1 /

\setsolid
\linethickness .2mm
\putrule from 9 -1 to 9 2.8
\putrule from 9 3.2 to  9 9
\putrule from 10 -1 to 10 3.8
\putrule from 10 4.2 to  10 9

\multiput{$\bullet$} at 6 6  7 5  8 4  10 4  /
\setsolid
\plot 6 6  8.8 3.2  / 
\plot 9.2 3.2  10 4 /

\put{$\circ$} at 9 3 

\multiput{} at 0 0  16 9 /
\put{(a)} at 8 -2.4

\endpicture} at 19 0
\put{\beginpicture
\setdots <1mm>
\plot -1 -1  9 9 /
\plot 9 8  17 -1 /

\setsolid
\linethickness .2mm
\putrule from 9 -1 to  9 0.8
\putrule from 9 1.2 to  9 9
\putrule from 10 0.2 to  10 9
\putrule from 10 -1 to  10 -.2

\multiput{$\bullet$} at 9 1  11 1  12 2  13 3  /
\setsolid
\plot 9 1  9.8 0.2  / 
\plot 10.2 0.2  13 3 /

\put{$\circ$} at 10 0

\multiput{} at 0 0  16 9 /
\put{(b)} at 8 -2.4

\endpicture} at 38 0
\put{\beginpicture
\setdots <1mm>
\plot -1 -1  9 9 /
\plot 9 8  17 -1 /

\setsolid
\linethickness .2mm
\putrule from 9 -1 to  9 9
\putrule from 10 -1 to  10 9

\multiput{$\bullet$} at 11 1  13 1  14 2   /
\setsolid
\plot 11 1  11.8  .2  / 
\plot 12.2 .2  14 2 /

\put{$\circ$} at 12 0 

\multiput{} at 0 0  16 9 /
\put{(b)} at 8 -2.4
\endpicture} at 57 0
\endpicture}
$$
	\medskip
Let $S$ be a subset of $\Bbb Z$. We say that a function $f\: S \to \Bbb Z$ 
is {\it eventually polynomial} of degree $t$ 
provided there is some integer $n_0$ such that
all integers $n \ge n_0$ belong to $S$ and the restriction of $f$
to the set $\{n\mid n\ge n_0\}$ is a polynomial function of degree $t$.
	\medskip
{\bf Corollary 2.} {\it The functions $d'_i, d''_i$ are eventually 
polynomial of degree $i$,
for any $i \ge 0$.} The corresponding polynomials 
are monic linear combinations of the binomial
coefficients $\binom t n$.
	\medskip
The proofs of Proposition 2 and Corollary 2 are similar to those of Proposition 1
and Corollary 1.
	\medskip
Here are the first polynomials which occur for $d'_i(t)$ and $d''_i(t)$:
$$
\align
 d'_0(t) &= 1 = \binom t0, \cr
 d'_1(t) &= t-1  = \binom t1 - \binom t0, \cr
 d'_2(t) &= \tfrac12(t^2-t-4)  = \binom t2 - 2 \binom t0, \cr
 d'_3(t) &= \tfrac 16(t^3-19t) = \binom t3 + \binom t2 - 3 \binom t1, \cr
  &\text{and} \cr
 d''_0(t) &= 1 = \binom t0, \cr
 d''_1(t) &= t = \binom t1, \cr
 d''_2(t) &= \tfrac12(t^2+t-8) = \binom t2 + \binom t1 - 4\binom t0,  \cr
 d''_3(t) &= \tfrac16(t^3+3t^2-28t-18) = \binom t3 + 2 \binom t2 - 4\binom t1 - 3 \binom                t0.
\endalign
$$
Again, we do not know a general formula for the polynomials $d'_i(t)$ and
$d''_i(t)$.

	\bigskip
Finally, the page with the odd-index Fibonacci partition triangle also
refers to some Fibonacci modules, namely the modules 
$R_t = R_t(x,y)$ considered in [FR2]. 
Again, it should be sufficient to discuss one example in detail,
say $R_4$. On the left, we exhibit the dimension vector of $R_4$, on the right
the corresponding two rows $t=3$ and $t=4$ of the ``triangle''. 
	
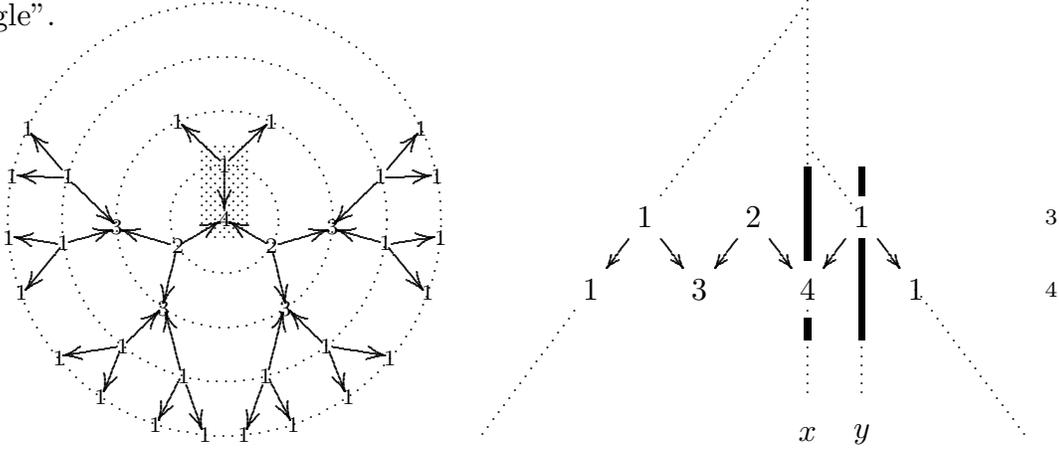
\bigskip
$$
{\beginpicture
\setcoordinatesystem units <0.6cm,0.6cm>

\multiput{\beginpicture

\put{$\ssize  4$} at 0 0
\multiput{$\ssize  2$} at -.866 -.5
  .866 -.5  /
\multiput{$\ssize  1$} at  0 1 /

\arrow <2mm> [0.25,0.75] from 0 0.9 to  0 0.2
\arrow <2mm> [0.25,0.75] from -.78 -.45
 to   -0.1 -0.05 
\arrow <2mm> [0.25,0.75] from .78 -.45
 to   0.1 -0.05 

\setshadegrid span <.5mm>
\vshade -.5 -.4 1.4  .5 -.4 1.4 /

\multiput{$\ssize  1$} at -.866 1.8 .866 1.8 /
\arrow <2mm> [0.25,0.75] from -0.1 1.1  to  -.78 1.75
\arrow <2mm> [0.25,0.75] from  0.1 1.1  to   .78 1.75




\startrotation by -0.5 0.866 about 0 0
\multiput{$\ssize  3$} at -.866 1.8 .866 1.8 /
\arrow <2mm> [0.25,0.75] from -0.1 1.1  to  -.78 1.75
\arrow <2mm> [0.25,0.75] from  0.1 1.1  to   .78 1.75
\multiput{$\ssize  1$} at -2.13 2.1  -1.1 2.8 
  2.13 2.1  1.1 2.8 /
\arrow <2mm> [0.25,0.75] from -2. 2.1     to  -.95 1.8 
\arrow <2mm> [0.25,0.75] from -1.07 2.7   to  -.866 1.9
\arrow <2mm> [0.25,0.75] from 2. 2.1      to  .95 1.8 
\arrow <2mm> [0.25,0.75] from 1.07 2.7    to  .866 1.9

\multiput{$\ssize  1$} at -3.27 2.3  -2.65 3 -1.7 3.62  -.7 3.93  /
\multiput{$\ssize  1$} at  3.27 2.3  2.65 3  1.7 3.62  .7 3.93 /

\arrow <2mm> [0.25,0.75] from -2.3 2.14   to  -3.15 2.3 
\arrow <2mm> [0.25,0.75] from -2.15 2.2    to  -2.6 2.9
\arrow <2mm> [0.25,0.75] from -1.15 2.9    to  -1.65 3.52
\arrow <2mm> [0.25,0.75] from -1.05 2.9    to  -.7 3.83 

\arrow <2mm> [0.25,0.75] from 2.3 2.14   to  3.15 2.3 
\arrow <2mm> [0.25,0.75] from 2.15 2.2    to  2.6 2.9
\arrow <2mm> [0.25,0.75] from 1.15 2.9    to  1.65 3.52
\arrow <2mm> [0.25,0.75] from 1.05 2.9    to  .7 3.83 
\stoprotation

\startrotation by -0.5 -0.866 about 0 0
\multiput{$\ssize  3$} at -.866 1.8 .866 1.8 /
\arrow <2mm> [0.25,0.75] from -0.1 1.1  to  -.78 1.75
\arrow <2mm> [0.25,0.75] from  0.1 1.1  to   .78 1.75
\multiput{$\ssize  1$} at -2.13 2.1  -1.1 2.8 
  2.13 2.1  1.1 2.8 /
\arrow <2mm> [0.25,0.75] from -2. 2.1     to  -.95 1.8 
\arrow <2mm> [0.25,0.75] from -1.07 2.7   to  -.866 1.9
\arrow <2mm> [0.25,0.75] from 2. 2.1      to  .95 1.8 
\arrow <2mm> [0.25,0.75] from 1.07 2.7    to  .866 1.9

\multiput{$\ssize  1$} at -3.27 2.3  -2.65 3 -1.7 3.62  -.7 3.93  /
\multiput{$\ssize  1$} at  3.27 2.3  2.65 3  1.7 3.62  .7 3.93 /

\arrow <2mm> [0.25,0.75] from -2.3 2.14   to  -3.15 2.3 
\arrow <2mm> [0.25,0.75] from -2.15 2.2    to  -2.6 2.9
\arrow <2mm> [0.25,0.75] from -1.15 2.9    to  -1.65 3.52
\arrow <2mm> [0.25,0.75] from -1.05 2.9    to  -.7 3.83 

\arrow <2mm> [0.25,0.75] from 2.3 2.14   to  3.15 2.3 
\arrow <2mm> [0.25,0.75] from 2.15 2.2    to  2.6 2.9
\arrow <2mm> [0.25,0.75] from 1.15 2.9    to  1.65 3.52
\arrow <2mm> [0.25,0.75] from 1.05 2.9    to  .7 3.83 
\stoprotation

\setdots <1mm>
\circulararc 360 degrees from 1 0 center at 0 0
\circulararc 360 degrees from 2 0 center at 0 0
\circulararc 360 degrees from 3 0 center at 0 0
\circulararc 360 degrees from 4 0 center at 0 0
\endpicture} at 0 0    /

\put{\beginpicture
\setcoordinatesystem units <.6cm,.8cm>
\put{$1$} at 0 0
\put{$1$} at 1 1
\put{$3$} at 2 0
\put{$2$} at 3 1
\put{$4$} at 4 0
\put{$1$} at 5 1
\put{$1$} at 6 0 
\arr{0.7 0.7}{0.3 0.3}
\arr{1.3 0.7}{1.7 0.3}
\arr{2.7 0.7}{2.3 0.3}
\arr{3.3 0.7}{3.7 0.3}
\arr{4.7 0.7}{4.3 0.3}
\arr{5.3 0.7}{5.7 0.3}

\linethickness.7mm
\putrule from 4 -.7 to 4 -.4
\putrule from 4  .4 to 4 1.7
\putrule from 5 -.7 to 5 .7
\putrule from 5  1.3 to 5 1.7

\put{$\ssize 3$} at 8.5 1
\put{$\ssize 4$} at 8.5 0

\setdots <1mm>
\plot 1.3 1.3  4 4  4 1.3 /
\plot -.3 -.3  -2 -2 /
\plot 4 -.3  4 -1.5 /
\plot 5 -.3  5 -1.5 /
\plot 4 2  5 1 /
\plot 6 0  8 -2 /
\put{$x$} at 4 -2
\put{$y$} at 5 -2
\endpicture} at 11 -.5

\endpicture}
$$
Whereas the left picture (copied from [FR2]) 
tries to put the emphasis on the vertex $x$ as the center,
the two pylons in the right picture stress that here we really deal with
a distribution of numbers which depend on an edge, thus on two neighboring
vertices, say $x$ and $y$. 
Note that the two pylons in the right picture mark the positions of
$x$ and $y$, the corresponding edge is indicated in the left picture by
the shaded region with the vertical arrow $y \to x$.
	\medskip
	
The numbers $d'_i(t)$ and $d''_i(t)$ 
are categorified by the modules $R_n = R_n(x,y)$; they
provide the Jordan-H\"older multiplicities of these modules.

First, let us consider the socle of $R_t$; the composition factors in
the socle of $R_t$ are of the form $S(z)$ with $0 \le D(x,z) \le t,$ 
$D(x,z) \equiv t \mod 2$, and $y \notin [x,z]$, or else with 
$0 < D(x,z) \le t-2,$ 
$D(x,z) \equiv t \mod 2$, and $y \in [x,z]$; here, $[x,z]$ denotes the 
path between $x$ and $z$. 
There is the following multiplicity formula:
	\medskip
\Rahmen
{$d'_i(t) = \dim R_t(x,y)_z \qquad \text{for}\quad D(x,z) = 
\left\{\matrix t-2i, &   & y\notin [x,z], \cr
                     &  \text{and} &      \cr
              2i-t, & & y\in [x,z], \endmatrix \right.$}
	\medskip
\noindent
We may reformulate the second case in terms of $d''_i(t)$ as follows:
	\medskip
\Rahmen
{$
 d''_i(t) = \dim R_t(x,y)_z \qquad \text{for}\quad D(x,z) = t-2i-2,\ y \in [x,z].
$}
	
	\medskip
(Remark: In the first displayed line, 
the condition $D(x,z) = t-2i$ implies that
we consider only 
$2i \le t$, the condition $D(x,z) = t-2i-2$ in the reformulation 
implies that $2i \le t-2$. Using the 
convention that $d'_i(t) = 0 = d''_i(t)$
for $i < 0$ and all $t$, we see again that these rules 
are valid also for $i < 0$.)
	\medskip
Note that the socle of $R_t(x,y)$ is isomorphic to the top of $R_{t+1}(x,y)$,
thus we can interprete $d'_i(t)$ and $d''_i(t)$ also as Jordan-H\"older multiplicities of composition factors in the top of  $R_{t+1}(x,y)$,

	\bigskip\bigskip
\vfill\eject
{\bf 3. Relations between the two triangles.}
	\medskip
Let us outline in which way the two triangles determine each other: 
in this way, we see that they are intimately connected.
	\medskip
{\bf Theorem.} For $t \ge 1$
$$
\alignat2
 d_i'(t) &=\ d_i(t)      - d_{i\!-\!1}(t\!-\!1)\qquad 
         &&\text{for all }i, \tag{1} \cr
 d_i''(t) &=\ d_{i\!+\!1}(t) - d_{i+1}(t\!-\!1)
         &&\text{for all }i,  \tag{2} \cr
 d_i(t\!-\!1) &=\ d'_{i+1}(t) - d'_{i+1}(t\!-\!1)  
         &&\text{for } i \le \tfrac{t-2}2,  \tag{3}\cr
 d_i(t\!-\!1) &=\ d''_i(t)   - d''_{i-1}(t\!-\!1) 
         &&\text{for } i \le \tfrac{t-2}2,  \tag{4} \cr
 d_i(t\!-\!1) &=\ d'_i(t)    - d''_{i-2}(t\!-\!1) 
         &&\text{for } i \le \tfrac{t}2.  \tag{5}
\endalignat
$$
	\medskip
Looking at the first four rules (1) to (4), 
one observes that the differences on the right
always concern differences along arrows: one starts with a suitable arrow, say $a \to b$ in one of the two triangles, and with 
the given additive function $g$ (either $g = d,$ or $g = d'$) and looks at the difference $g(b)-g(a).$ In contrast, the rule (5) deals with
differences of numbers often quite far apart; of special interest are the
cases $t = 2i$ where one deals with the difference along a knight's move
(known from chess):
$$
 d_i(2i) = d'_i(2i\!+\!1)\ -d''_{i-2}(2i).
$$
Say for $i= 3$, one considers 
$$
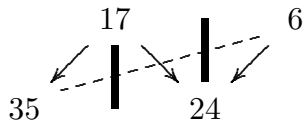
{\beginpicture
\setcoordinatesystem units <1cm,1cm>
\put{$35$} at 0 0
\put{$17$} at 1 1
\put{$24$} at 2 0
\put{$6$} at 3 1
\arr{0.7 0.7}{0.3 0.3}
\arr{1.3 0.7}{1.7 0.3}
\arr{2.7 0.7}{2.3 0.3}
\linethickness.7mm
\putrule from 1 0  to 1 0.7
\putrule from 2 0.3  to 2 1
\setdashes <1mm>
\plot 2.6 0.8  0.4 0.2 /
\endpicture}
$$
here, $29 = d_3(6) = d'_3(7) - d''_1(6) = 35-6.$ 
These differences along knight's moves yield all  
the numbers on the pylon of the even-index triangle.

Starting with the even-index triangle, (1) and (2) assert that
we obtain the numbers $d'_i$ as differences along the 
south-west arrows, the numbers $d''_i$ as differences along the
south-east arrows.
Similarly, the rules (3), (4) and (5) show how to obtain the numbers $d_i(t)$ 
of the even-index triangle from the odd-index triangle. 
Those outside of the pylon
are obtained in three different ways, twice as differences 
along the south-west arrows, namely, on the one hand,
looking at the arrows left of the second pylon, see (3), and, on the other hand,
also by looking at the arrows right of the first pylon, see (4). But
all the numbers of the even-index triangle are obtained using the rule (5).
	\medskip\medskip
\vfill\eject
{\bf Differences on arrows for the even-index triangle.}

For the convenience of the reader, let us visualize the rules (1) and (2)
by exhibiting part of the even-index triangle and inserting the 
difference numbers to some of the arrows:
$$
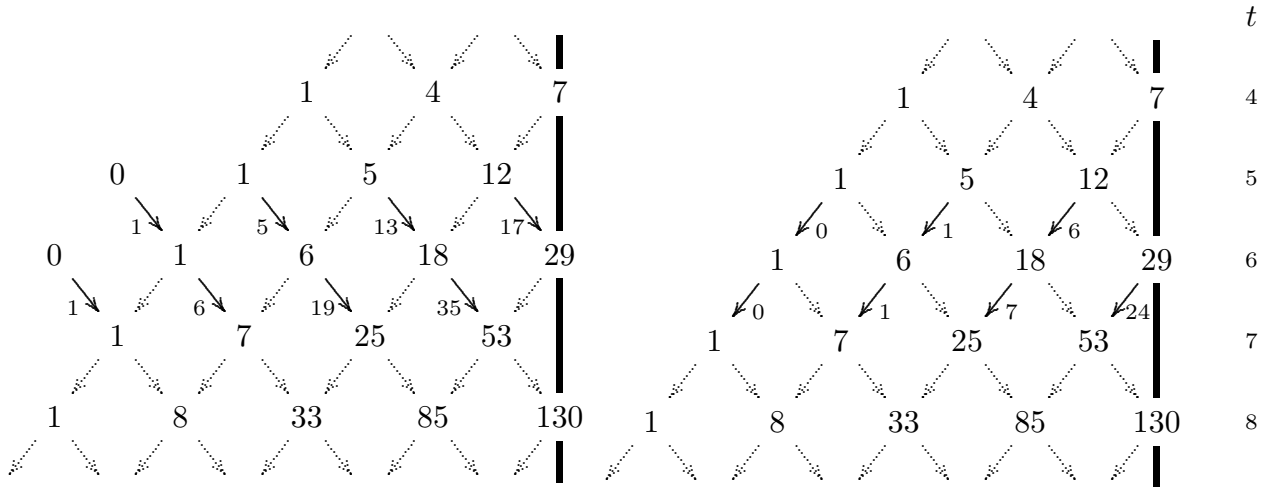
{\beginpicture
\setcoordinatesystem units <.7cm,.9cm>

\put{\beginpicture
\multiput{} at -6 -3  1.5 -9 /

\multiput{$1$} at   -4 -4  -5 -5  -6 -6  -7 -7  -8 -8  /

\put{$0$} at -7 -5 
\put{$0$} at -8 -6

\put{$4$} at -2 -4 
\put{$5$} at -3 -5 
\put{$6$} at -4 -6 
\put{$7$} at -5 -7 
\put{$8$} at -6 -8 

\put{$7$} at 0 -4
\put{$12$} at -1 -5
\put{$18$} at -2 -6 
\put{$25$} at -3 -7 
\put{$33$} at -4 -8 

\put{$29$} at 0 -6
\put{$53$} at -1 -7 
\put{$85$} at -2 -8 

\put{$130$} at 0 -8

\setdots <.5mm>

\arr{-3.3 -3.3}{-3.7 -3.7}
\arr{-2.7 -3.3}{-2.3 -3.7}
\arr{-1.3 -3.3}{-1.7 -3.7}
\arr{-.7 -3.3}{-.3 -3.7}

\arr{-4.3 -4.3}{-4.7 -4.7}
\arr{-3.7 -4.3}{-3.3 -4.7}
\arr{-2.3 -4.3}{-2.7 -4.7}
\arr{-1.7 -4.3}{-1.3 -4.7}
\arr{-.3 -4.3}{-.7 -4.7}

\arr{-5.3 -5.3}{-5.7 -5.7}
\arr{-3.3 -5.3}{-3.7 -5.7}
\arr{-1.3 -5.3}{-1.7 -5.7}

\arr{-6.3 -6.3}{-6.7 -6.7}
\arr{-4.3 -6.3}{-4.7 -6.7}
\arr{-2.3 -6.3}{-2.7 -6.7}
\arr{-.3 -6.3}{-.7 -6.7}

\setsolid

\arr{-6.7 -5.3}{-6.3 -5.7}
\arr{-4.7 -5.3}{-4.3 -5.7}
\arr{-2.7 -5.3}{-2.3 -5.7}
\arr{-.7 -5.3}{-.3 -5.7}

\arr{-7.7 -6.3}{-7.3 -6.7}
\arr{-5.7 -6.3}{-5.3 -6.7}
\arr{-3.7 -6.3}{-3.3 -6.7}
\arr{-1.7 -6.3}{-1.3 -6.7}

\setdots<.5mm>

\arr{-7.3 -7.3}{-7.7 -7.7}
\arr{-6.7 -7.3}{-6.3 -7.7}
\arr{-5.3 -7.3}{-5.7 -7.7}
\arr{-4.7 -7.3}{-4.3 -7.7}
\arr{-3.3 -7.3}{-3.7 -7.7}
\arr{-2.7 -7.3}{-2.3 -7.7}
\arr{-1.3 -7.3}{-1.7 -7.7}
\arr{-.7 -7.3}{-.3 -7.7}

\arr{-8.3 -8.3}{-8.7 -8.7}
\arr{-7.7 -8.3}{-7.3 -8.7}
\arr{-6.3 -8.3}{-6.7 -8.7}
\arr{-5.7 -8.3}{-5.3 -8.7}
\arr{-4.3 -8.3}{-4.7 -8.7}
\arr{-3.7 -8.3}{-3.3 -8.7}
\arr{-2.3 -8.3}{-2.7 -8.7}
\arr{-1.7 -8.3}{-1.3 -8.7}
\arr{-.3 -8.3}{-.7 -8.7}
\setsolid

\linethickness.7mm
\putrule from 0 -3.3 to 0 -3.7
\putrule from 0 -4.3 to 0 -5.7
\putrule from 0 -6.3 to 0 -7.7
\putrule from 0 -8.3 to 0 -8.8

\put{$\ssize 1$}  at -6.7 -5.65
\put{$\ssize 5$}  at -4.7 -5.65
\put{$\ssize 13$} at -2.75 -5.65
\put{$\ssize 17$} at -.75 -5.65

\put{$\ssize 1$}  at -7.7 -6.65
\put{$\ssize 6$}  at -5.7 -6.65
\put{$\ssize 19$} at -3.75 -6.65
\put{$\ssize 35$} at -1.75 -6.65

\endpicture} at 0 0
\put{\beginpicture
\multiput{} at -6 -3  1.5 -9 /

\multiput{$1$} at   -4 -4  -5 -5  -6 -6  -7 -7  -8 -8  /

\put{$4$} at -2 -4 
\put{$5$} at -3 -5 
\put{$6$} at -4 -6 
\put{$7$} at -5 -7 
\put{$8$} at -6 -8 

\put{$7$} at 0 -4
\put{$12$} at -1 -5
\put{$18$} at -2 -6 
\put{$25$} at -3 -7 
\put{$33$} at -4 -8 

\put{$29$} at 0 -6
\put{$53$} at -1 -7 
\put{$85$} at -2 -8 

\put{$130$} at 0 -8

\setdots <.5mm>

\arr{-3.3 -3.3}{-3.7 -3.7}
\arr{-2.7 -3.3}{-2.3 -3.7}
\arr{-1.3 -3.3}{-1.7 -3.7}
\arr{-.7 -3.3}{-.3 -3.7}

\arr{-4.3 -4.3}{-4.7 -4.7}
\arr{-3.7 -4.3}{-3.3 -4.7}
\arr{-2.3 -4.3}{-2.7 -4.7}
\arr{-1.7 -4.3}{-1.3 -4.7}
\arr{-.3 -4.3}{-.7 -4.7}

\arr{-4.7 -5.3}{-4.3 -5.7}
\arr{-2.7 -5.3}{-2.3 -5.7}
\arr{-.7 -5.3}{-.3 -5.7}

\arr{-5.7 -6.3}{-5.3 -6.7}
\arr{-3.7 -6.3}{-3.3 -6.7}
\arr{-1.7 -6.3}{-1.3 -6.7}

\setsolid

\arr{-5.3 -5.3}{-5.7 -5.7}
\arr{-3.3 -5.3}{-3.7 -5.7}
\arr{-1.3 -5.3}{-1.7 -5.7}

\arr{-6.3 -6.3}{-6.7 -6.7}
\arr{-4.3 -6.3}{-4.7 -6.7}
\arr{-2.3 -6.3}{-2.7 -6.7}
\arr{-.3 -6.3}{-.7 -6.7}

\setdots<.5mm>

\arr{-7.3 -7.3}{-7.7 -7.7}
\arr{-6.7 -7.3}{-6.3 -7.7}
\arr{-5.3 -7.3}{-5.7 -7.7}
\arr{-4.7 -7.3}{-4.3 -7.7}
\arr{-3.3 -7.3}{-3.7 -7.7}
\arr{-2.7 -7.3}{-2.3 -7.7}
\arr{-1.3 -7.3}{-1.7 -7.7}
\arr{-.7 -7.3}{-.3 -7.7}

\arr{-8.3 -8.3}{-8.7 -8.7}
\arr{-7.7 -8.3}{-7.3 -8.7}
\arr{-6.3 -8.3}{-6.7 -8.7}
\arr{-5.7 -8.3}{-5.3 -8.7}
\arr{-4.3 -8.3}{-4.7 -8.7}
\arr{-3.7 -8.3}{-3.3 -8.7}
\arr{-2.3 -8.3}{-2.7 -8.7}
\arr{-1.7 -8.3}{-1.3 -8.7}
\arr{-.3 -8.3}{-.7 -8.7}
\setsolid

\linethickness.7mm
\putrule from 0 -3.3 to 0 -3.7
\putrule from 0 -4.3 to 0 -5.7
\putrule from 0 -6.3 to 0 -7.7
\putrule from 0 -8.3 to 0 -8.8

\put{$t$} at 1.5 -3
\put{$\ssize 4$} at 1.5 -4
\put{$\ssize 5$} at 1.5 -5
\put{$\ssize 6$} at 1.5 -6
\put{$\ssize 7$} at 1.5 -7
\put{$\ssize 8$} at 1.5 -8

\put{$\ssize 0$}  at -5.3 -5.65
\put{$\ssize 1$}  at -3.3 -5.65
\put{$\ssize 6$}  at -1.3 -5.65

\put{$\ssize 0$}  at -6.3 -6.65
\put{$\ssize 1$}  at -4.3 -6.65
\put{$\ssize 7$}  at -2.3 -6.65
\put{$\ssize 24$} at -.3 -6.65

\endpicture} at 9.5 0

\endpicture}
$$

In a slightly different way, we may start with $d$ as
an additive function on $\Bbb N\Delta^{\odd}$ and consider only the differences
along the south-east arrows. (It is quite important to be aware that
any additive function $g$ on $\Bbb N\Delta^{\ev}$ gives rise to an additive
function on $\Bbb N\Delta^{\odd}$, also denoted by $g$, by extending the given
function via the rule $g(i,t) = g(t-i,t).$)
In the following display 
we have inserted the differences for all the south-east arrows ending in the
layers $t=0,1,2,3$ and $t=6,7.$  
$$
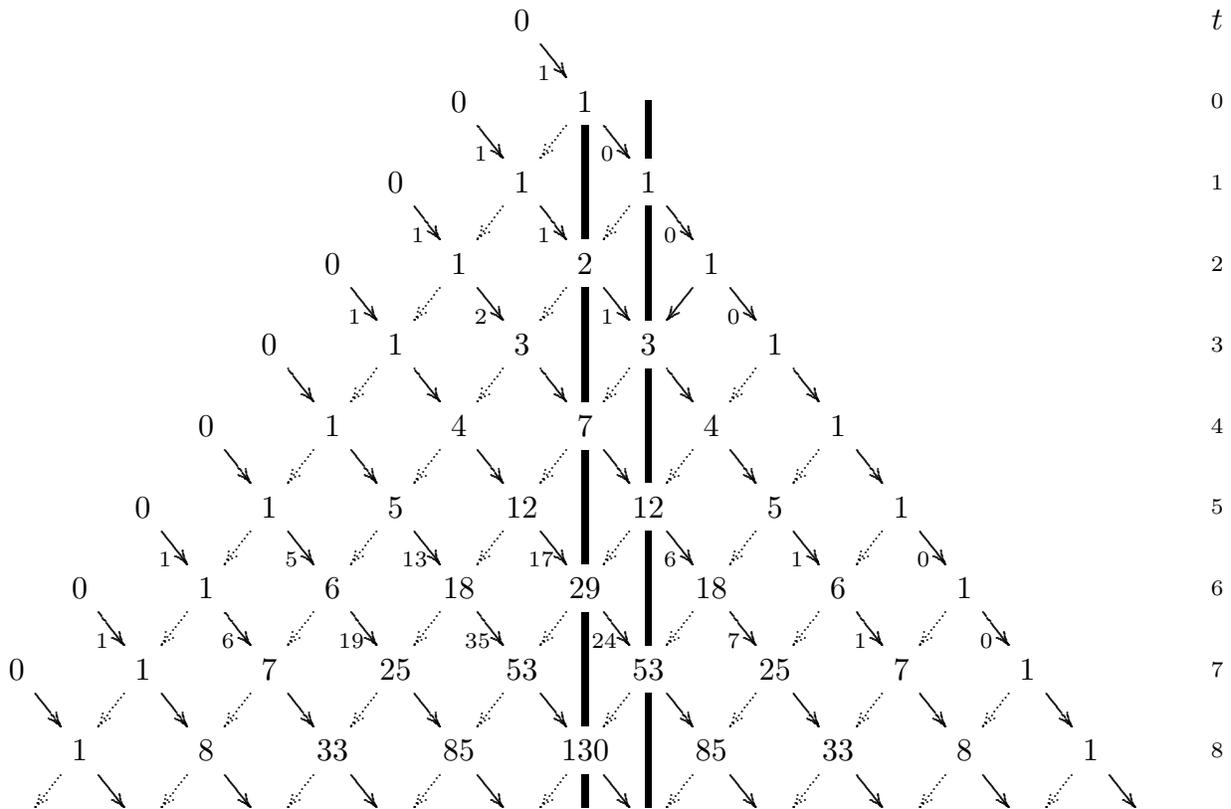
{\beginpicture
\setcoordinatesystem units <.7cm,.9cm>

\multiput{$1$} at    -5 -5  -6 -6  -7 -7  -8 -8  -9 -9 /

\put{$4$} at -3 -5 
\put{$5$} at -4 -6 
\put{$6$} at -5 -7 
\put{$7$} at -6 -8 
\put{$8$} at -7 -9 

\put{$7$} at -1 -5
\put{$12$} at -2 -6 
\put{$18$} at -3 -7 
\put{$25$} at -4 -8 
\put{$33$} at -5 -9 

\put{$12$} at 0 -6
\put{$29$} at -1 -7 
\put{$53$} at -2 -8 
\put{$85$} at -3 -9 

\put{$130$} at -1 -9

\multiput{$1$} at -4 -4  -3 -3 -2 -2  -1 -1  0 -2  1 -3  2 -4 /
\multiput{$2$} at -1 -3 /
\multiput{$3$} at -2 -4  0 -4 /

\setsolid 

\arr{-1.7 -.3}{-1.3 -.7}

\arr{-2.7 -1.3}{-2.3 -1.7}
\arr{-.7 -1.3}{-.3 -1.7}

\arr{-3.7 -2.3}{-3.3 -2.7}
\arr{-1.7 -2.3}{-1.3 -2.7}

\arr{-4.7 -3.3}{-4.3 -3.7}
\arr{-2.7 -3.3}{-2.3 -3.7}
\arr{-.7 -3.3}{-.3 -3.7}
\arr{.7 -3.3}{.3 -3.7}

\arr{-3.7 -4.3}{-3.3 -4.7}
\arr{-1.7 -4.3}{-1.3 -4.7}

\arr{.3 -2.3}{.7 -2.7}

\arr{1.3 -3.3}{1.7 -3.7}

\setdots <.5mm> 
\arr{-1.3 -1.3}{-1.7 -1.7}
\arr{-2.3 -2.3}{-2.7 -2.7}
\arr{-.3 -2.3}{-.7 -2.7}

\arr{-3.3 -3.3}{-3.7 -3.7}
\arr{-1.3 -3.3}{-1.7 -3.7}
\arr{-2.3 -4.3}{-2.7 -4.7}
\arr{-.3 -4.3}{-.7 -4.7}
\arr{-4.3 -4.3}{-4.7 -4.7}

\arr{-6.3 -6.3}{-6.7 -6.7}
\arr{-4.3 -6.3}{-4.7 -6.7}
\arr{-2.3 -6.3}{-2.7 -6.7}
\arr{-.3 -6.3}{-.7 -6.7}

\arr{-7.3 -7.3}{-7.7 -7.7}
\arr{-5.3 -7.3}{-5.7 -7.7}
\arr{-3.3 -7.3}{-3.7 -7.7}
\arr{-1.3 -7.3}{-1.7 -7.7}

\arr{-5.3 -5.3}{-5.7 -5.7}
\arr{-3.3 -5.3}{-3.7 -5.7}
\arr{-1.3 -5.3}{-1.7 -5.7}

\arr{-8.3 -8.3}{-8.7 -8.7}
\arr{-6.3 -8.3}{-6.7 -8.7}
\arr{-4.3 -8.3}{-4.7 -8.7}
\arr{-2.3 -8.3}{-2.7 -8.7}
\arr{-.3 -8.3}{-.7 -8.7}
\arr{-9.3 -9.3}{-9.7 -9.7}
\arr{-7.3 -9.3}{-7.7 -9.7}
\arr{-5.3 -9.3}{-5.7 -9.7}
\arr{-3.3 -9.3}{-3.7 -9.7}
\arr{-1.3 -9.3}{-1.7 -9.7}

\arr{1.7 -6.3}{1.3 -6.7}
\arr{3.7 -6.3}{3.3 -6.7}

\arr{4.7 -7.3}{4.3 -7.7}
\arr{2.7 -7.3}{2.3 -7.7}
\arr{.7 -7.3}{.3 -7.7}
\arr{1.7 -4.3}{1.3 -4.7}
\arr{2.7 -5.3}{2.3 -5.7}
\arr{.7 -5.3}{.3 -5.7}
\arr{1.7 -8.3}{1.3 -8.7}
\arr{5.7 -8.3}{5.3 -8.7}
\arr{3.7 -8.3}{3.3 -8.7}
\arr{6.7 -9.3}{6.3 -9.7}
\arr{4.7 -9.3}{4.3 -9.7}
\arr{2.7 -9.3}{2.3 -9.7}
\arr{.7 -9.3}{.3 -9.7}

\setsolid
\arr{-4.7 -5.3}{-4.3 -5.7}
\arr{-2.7 -5.3}{-2.3 -5.7}
\arr{-.7 -5.3}{-.3 -5.7}
\arr{-5.7 -6.3}{-5.3 -6.7}
\arr{-3.7 -6.3}{-3.3 -6.7}
\arr{-1.7 -6.3}{-1.3 -6.7}
\arr{-6.7 -7.3}{-6.3 -7.7}
\arr{-4.7 -7.3}{-4.3 -7.7}
\arr{-2.7 -7.3}{-2.3 -7.7}
\arr{-.7 -7.3}{-.3 -7.7}

\arr{-7.7 -8.3}{-7.3 -8.7}
\arr{-5.7 -8.3}{-5.3 -8.7}
\arr{-3.7 -8.3}{-3.3 -8.7}

\arr{-1.7 -8.3}{-1.3 -8.7}

\arr{-8.7 -9.3}{-8.3 -9.7}
\arr{-6.7 -9.3}{-6.3 -9.7}
\arr{-4.7 -9.3}{-4.3 -9.7}
\arr{-2.7 -9.3}{-2.3 -9.7}
\arr{-.7 -9.3}{-.3 -9.7}

\setsolid

\linethickness.7mm
\putrule from 0 -1. to 0 -1.7
\putrule from 0 -2.3 to 0 -3.7
\putrule from 0 -4.3 to 0 -5.7
\putrule from 0 -6.3 to 0 -7.7
\putrule from 0 -8.3 to 0 -9.7

\putrule from -1 -1.3 to -1 -2.7
\putrule from -1 -3.3 to -1 -4.7
\putrule from -1 -5.3 to -1 -6.7
\putrule from -1 -7.3 to -1 -8.7
\putrule from -1 -9.3 to -1 -9.7

\put{$\ssize 0$} at 9 -1
\put{$\ssize 1$} at 9 -2
\put{$\ssize 2$} at 9 -3
\put{$\ssize 3$} at 9 -4
\put{$\ssize 4$} at 9 -5
\put{$\ssize 5$} at 9 -6
\put{$\ssize 6$} at 9 -7
\put{$\ssize 7$} at 9 -8
\put{$\ssize 8$} at 9 -9

\put{$t$} at 9 0

\arr{.3 -4.3}{.7 -4.7}
\arr{2.3 -4.3}{2.7 -4.7}

\arr{3.3 -5.3}{3.7 -5.7}

\arr{1.3 -5.3}{1.7 -5.7}
\arr{2.3 -6.3}{2.7 -6.7}
\arr{.3 -6.3}{.7 -6.7}
\arr{4.3 -6.3}{4.7 -6.7}

\arr{5.3 -7.3}{5.7 -7.7}
\arr{3.3 -7.3}{3.7 -7.7}
\arr{1.3 -7.3}{1.7 -7.7}

\arr{6.3 -8.3}{6.7 -8.7}
\arr{4.3 -8.3}{4.7 -8.7}
\arr{2.3 -8.3}{2.7 -8.7}
\arr{.3 -8.3}{.7 -8.7}

\arr{7.3 -9.3}{7.7 -9.7}
\arr{5.3 -9.3}{5.7 -9.7}
\arr{3.3 -9.3}{3.7 -9.7}
\arr{1.3 -9.3}{1.7 -9.7}

\put{$4$} at 1 -5  
\put{$5$} at 2 -6  
\put{$6$} at 3 -7  
\put{$7$} at 4 -8  
\put{$8$} at 5 -9

\put{$18$} at 1 -7  
\put{$25$} at 2 -8
\put{$33$} at 3 -9

\put{$53$} at 0 -8  
\put{$85$} at 1 -9

\multiput{$1$} at 3 -5  4 -6  5 -7  6 -8  7 -9 /

\multiput{$0$} at -10 -8  
  -9 -7  -8 -6  -7 -5  -6 -4  -5 -3  -4 -2  -3 -1  -2 0 /

\arr{-5.7 -4.3}{-5.3 -4.7}
\arr{-6.7 -5.3}{-6.3 -5.7}

\arr{-7.7 -6.3}{-7.3 -6.7}
\arr{-8.7 -7.3}{-8.3 -7.7}
\arr{-9.7 -8.3}{-9.3 -8.7}

\put{$\ssize 1$} at -1.65  -.65

\put{$\ssize 1$} at -2.65  -1.65
\put{$\ssize 0$} at  -.65  -1.65

\put{$\ssize 1$} at -3.65  -2.65
\put{$\ssize 1$} at -1.65  -2.65
\put{$\ssize 0$} at  .35  -2.65

\put{$\ssize 1$} at -4.65  -3.65
\put{$\ssize 2$} at -2.65  -3.65
\put{$\ssize 1$} at  -.65  -3.65
\put{$\ssize 0$} at  1.35  -3.65

\put{$\ssize 1$} at -7.65  -6.65
\put{$\ssize 5$} at -5.65  -6.65
\put{$\ssize 13$} at -3.7  -6.65
\put{$\ssize 17$} at -1.7  -6.65
\put{$\ssize 6$} at  .35  -6.65
\put{$\ssize 1$} at  2.35  -6.65
\put{$\ssize 0$} at  4.35  -6.65

\put{$\ssize 1$} at -8.65  -7.65
\put{$\ssize 6$} at -6.65  -7.65
\put{$\ssize 19$} at -4.7  -7.65
\put{$\ssize 35$} at -2.7  -7.65
\put{$\ssize 24$} at  -.7  -7.65
\put{$\ssize 7$} at  1.35  -7.65
\put{$\ssize 1$} at  3.35  -7.65
\put{$\ssize 0$} at  5.35  -7.65

\comment
\setdots <1mm>
\plot -12.5 -1  -4 -1 /
\plot -12.5 -2  -5 -2 /
\plot -12.5 -3  -6 -3 /
\plot -12.5 -4  -7 -4 /
\plot -12.5 -6  -9 -6 /
\plot -12.5 -7  -10 -7 /
\plot -12.5 -8  -10 -8 /

\put{$R_6$} at -12 -6.5
\put{$R_7$} at -12 -7.5

\put{$R_1$} at -12 -1.5
\put{$R_2$} at -12 -2.5
\put{$R_3$} at -12 -3.5
\endcomment

\endpicture}
$$
What do we see? Looking at the difference numbers, we obtain precisely
the odd-index triangle. The only deviation may concern 
the position of the pylons;
they should be shifted slightly to the left; but we can 
overcome this problem
by inserting the difference number for any arrow
$a \to b$ directly at the position $b$; then also the position of the pylons
is correct.

\bigskip

{\bf Differences for the odd-index triangle.}

In the even-index pictures, we have used differences along the
given arrows (as well as additional ones on the left
boundary), according to the rules (1) and (2). Now we start with
the odd-index triangle, looking at the south-west arrows (as
well as the corresponding arrows on the right boundary), see
the rules (3) and (4):
$$
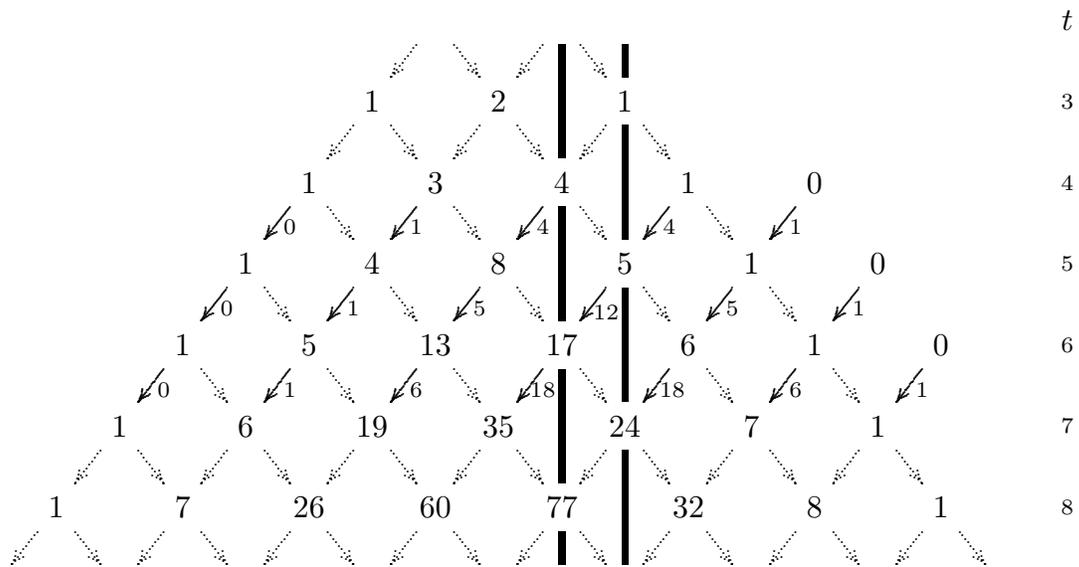
{\beginpicture
\setcoordinatesystem units <.7cm,.9cm>

\multiput{$1$} at  -4 -4  -5 -5  -6 -6  -7 -7  -8 -8  -9 -9 /

\put{$2$} at -2 -4 
\put{$3$} at -3 -5 
\put{$4$} at -4 -6 
\put{$5$} at -5 -7 
\put{$6$} at -6 -8 
\put{$7$} at -7 -9 

\put{$1$} at 0 -4
\put{$4$} at -1 -5
\put{$8$} at -2 -6 
\put{$13$} at -3 -7 
\put{$19$} at -4 -8 
\put{$26$} at -5 -9 

\put{$5$} at 0 -6
\put{$17$} at -1 -7 
\put{$35$} at -2 -8 
\put{$60$} at -3 -9 

\put{$77$} at -1 -9 

\setdots <.5mm>

\arr{-3.3 -3.3}{-3.7 -3.7}
\arr{-2.7 -3.3}{-2.3 -3.7}
\arr{-1.3 -3.3}{-1.7 -3.7}
\arr{-.7 -3.3}{-.3 -3.7}

\arr{-4.3 -4.3}{-4.7 -4.7}
\arr{-3.7 -4.3}{-3.3 -4.7}
\arr{-2.3 -4.3}{-2.7 -4.7}
\arr{-1.7 -4.3}{-1.3 -4.7}
\arr{-.3 -4.3}{-.7 -4.7}

\setsolid

\arr{-5.3 -5.3}{-5.7 -5.7}
\arr{-3.3 -5.3}{-3.7 -5.7}
\arr{-1.3 -5.3}{-1.7 -5.7}

\arr{-6.3 -6.3}{-6.7 -6.7}
\arr{-4.3 -6.3}{-4.7 -6.7}
\arr{-2.3 -6.3}{-2.7 -6.7}
\arr{-.3 -6.3}{-.7 -6.7}

\arr{-7.3 -7.3}{-7.7 -7.7}
\arr{-5.3 -7.3}{-5.7 -7.7}
\arr{-3.3 -7.3}{-3.7 -7.7}
\arr{-1.3 -7.3}{-1.7 -7.7}

\setdots <.5mm>

\arr{-4.7 -5.3}{-4.3 -5.7}
\arr{-2.7 -5.3}{-2.3 -5.7}
\arr{-.7 -5.3}{-.3 -5.7}
\arr{-5.7 -6.3}{-5.3 -6.7}
\arr{-3.7 -6.3}{-3.3 -6.7}
\arr{-1.7 -6.3}{-1.3 -6.7}
\arr{-6.7 -7.3}{-6.3 -7.7}
\arr{-4.7 -7.3}{-4.3 -7.7}
\arr{-2.7 -7.3}{-2.3 -7.7}
\arr{-.7 -7.3}{-.3 -7.7}

\arr{-8.3 -8.3}{-8.7 -8.7}
\arr{-7.7 -8.3}{-7.3 -8.7}
\arr{-6.3 -8.3}{-6.7 -8.7}
\arr{-5.7 -8.3}{-5.3 -8.7}
\arr{-4.3 -8.3}{-4.7 -8.7}
\arr{-3.7 -8.3}{-3.3 -8.7}
\arr{-2.3 -8.3}{-2.7 -8.7}
\arr{-1.7 -8.3}{-1.3 -8.7}
\arr{-.3 -8.3}{-.7 -8.7}

\arr{-9.3 -9.3}{-9.7 -9.7}
\arr{-8.7 -9.3}{-8.3 -9.7}
\arr{-7.3 -9.3}{-7.7 -9.7}
\arr{-6.7 -9.3}{-6.3 -9.7}
\arr{-5.3 -9.3}{-5.7 -9.7}
\arr{-4.7 -9.3}{-4.3 -9.7}
\arr{-3.3 -9.3}{-3.7 -9.7}
\arr{-2.7 -9.3}{-2.3 -9.7}
\arr{-1.3 -9.3}{-1.7 -9.7}
\arr{-.7 -9.3}{-.3 -9.7}

\setsolid

\linethickness.7mm
\putrule from 0 -3.3 to 0 -3.7
\putrule from 0 -4.3 to 0 -5.7
\putrule from 0 -6.3 to 0 -7.7
\putrule from 0 -8.3 to 0 -9.7

\putrule from -1 -3.3 to -1 -4.7
\putrule from -1 -5.3 to -1 -6.7
\putrule from -1 -7.3 to -1 -8.7
\putrule from -1 -9.3 to -1 -9.7

\put{$\ssize 3$} at 7 -4
\put{$\ssize 4$} at 7 -5
\put{$\ssize 5$} at 7 -6
\put{$\ssize 6$} at 7 -7
\put{$\ssize 7$} at 7 -8
\put{$\ssize 8$} at 7 -9

\put{$t$} at 7 -3

\arr{2.7 -5.3}{2.3 -5.7}
\arr{.7 -5.3}{.3 -5.7}

\arr{3.7 -6.3}{3.3 -6.7}
\arr{1.7 -6.3}{1.3 -6.7}

\arr{4.7 -7.3}{4.3 -7.7}
\arr{2.7 -7.3}{2.3 -7.7}
\arr{.7 -7.3}{.3 -7.7}

\setdots <.5mm>
\arr{.3 -4.3}{.7 -4.7}

\arr{1.3 -5.3}{1.7 -5.7}
\arr{2.3 -6.3}{2.7 -6.7}
\arr{.3 -6.3}{.7 -6.7}
\arr{3.3 -7.3}{3.7 -7.7}
\arr{1.3 -7.3}{1.7 -7.7}

\arr{4.3 -8.3}{4.7 -8.7}
\arr{3.7 -8.3}{3.3 -8.7}
\arr{2.3 -8.3}{2.7 -8.7}
\arr{1.7 -8.3}{1.3 -8.7}
\arr{.3 -8.3}{.7 -8.7}

\arr{5.3 -9.3}{5.7 -9.7}
\arr{4.7 -9.3}{4.3 -9.7}
\arr{3.3 -9.3}{3.7 -9.7}
\arr{2.7 -9.3}{2.3 -9.7}
\arr{1.3 -9.3}{1.7 -9.7}
\arr{.7 -9.3}{.3 -9.7}

\put{$1$} at 1 -5  
\put{$1$} at 2 -6  
\put{$1$} at 3 -7  
\put{$1$} at 4 -8  
\put{$1$} at 5 -9

\put{$6$} at 1 -7  
\put{$7$} at 2 -8
\put{$8$} at 3 -9

\put{$24$} at 0 -8  
\put{$32$} at 1 -9  

\comment
\setdots <1mm>
\plot -10.5 -5  -6 -5 /
\plot -10.5 -6  -7 -6 /
\plot -10.5 -7  -8 -7 /
\plot -10.5 -8  -9 -8 /

\put{$R_5$} at -10 -5.5
\put{$R_6$} at -10 -6.5
\put{$R_7$} at -10 -7.5
\endcomment

\multiput{$0$} at 3 -5  4 -6  5 -7 /

\put{$\ssize 0$} at -5.3  -5.55
\put{$\ssize 1$} at -3.3  -5.55
\put{$\ssize 4$} at -1.3  -5.55
\put{$\ssize 4$} at   .7  -5.55
\put{$\ssize 1$} at  2.7  -5.55

\put{$\ssize 0$} at -6.3  -6.55
\put{$\ssize 1$} at -4.3  -6.55
\put{$\ssize 5$} at -2.3  -6.55
\put{$\ssize 12$} at -.3  -6.6
\put{$\ssize 5$} at  1.7  -6.55
\put{$\ssize 1$} at  3.7  -6.55

\put{$\ssize 0$} at -7.3  -7.55
\put{$\ssize 1$} at -5.3  -7.55
\put{$\ssize 6$} at -3.3  -7.55
\put{$\ssize 18$} at -1.3  -7.55
\put{$\ssize 18$} at   .75  -7.55
\put{$\ssize 6$} at  2.7  -7.55
\put{$\ssize 1$} at  4.7  -7.55

\endpicture}
$$
	\medskip
The rule (5) deals with the following differences:
$$
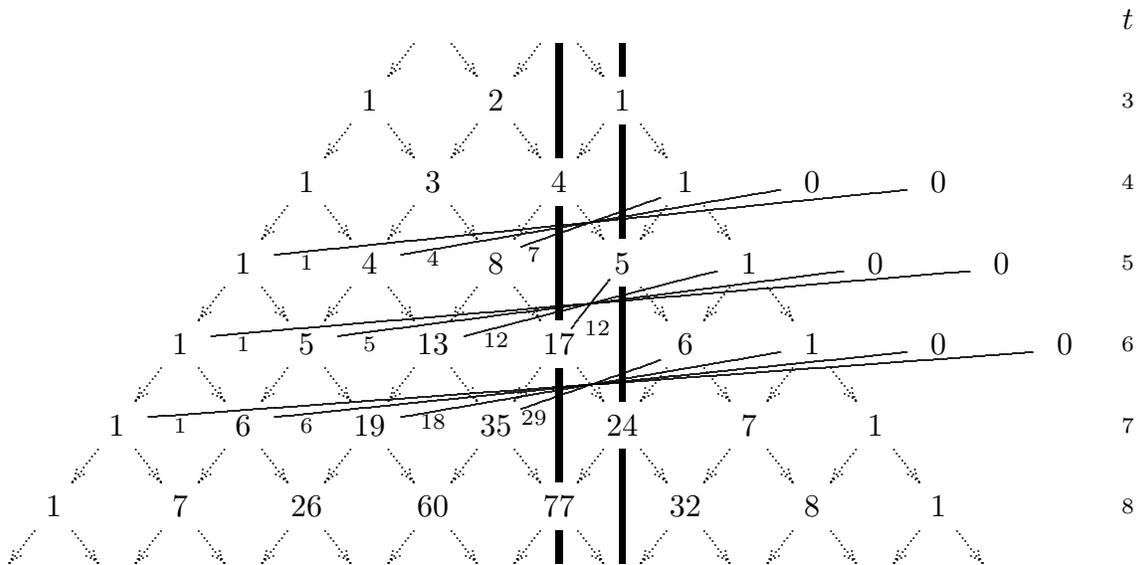
{\beginpicture
\setcoordinatesystem units <.7cm,.9cm>

\multiput{$1$} at  -4 -4  -5 -5  -6 -6  -7 -7  -8 -8  -9 -9 /

\put{$2$} at -2 -4 
\put{$3$} at -3 -5 
\put{$4$} at -4 -6 
\put{$5$} at -5 -7 
\put{$6$} at -6 -8 
\put{$7$} at -7 -9 

\put{$1$} at 0 -4
\put{$4$} at -1 -5
\put{$8$} at -2 -6 
\put{$13$} at -3 -7 
\put{$19$} at -4 -8 
\put{$26$} at -5 -9 

\put{$5$} at 0 -6
\put{$17$} at -1 -7 
\put{$35$} at -2 -8 
\put{$60$} at -3 -9 

\put{$77$} at -1 -9 

\setdots <.5mm>

\arr{-3.3 -3.3}{-3.7 -3.7}
\arr{-2.7 -3.3}{-2.3 -3.7}
\arr{-1.3 -3.3}{-1.7 -3.7}
\arr{-.7 -3.3}{-.3 -3.7}

\arr{-4.3 -4.3}{-4.7 -4.7}
\arr{-3.7 -4.3}{-3.3 -4.7}
\arr{-2.3 -4.3}{-2.7 -4.7}
\arr{-1.7 -4.3}{-1.3 -4.7}
\arr{-.3 -4.3}{-.7 -4.7}

\arr{-5.3 -5.3}{-5.7 -5.7}
\arr{-3.3 -5.3}{-3.7 -5.7}
\arr{-1.3 -5.3}{-1.7 -5.7}

\arr{-6.3 -6.3}{-6.7 -6.7}
\arr{-4.3 -6.3}{-4.7 -6.7}
\arr{-2.3 -6.3}{-2.7 -6.7}

\arr{-7.3 -7.3}{-7.7 -7.7}
\arr{-5.3 -7.3}{-5.7 -7.7}
\arr{-3.3 -7.3}{-3.7 -7.7}
\arr{-1.3 -7.3}{-1.7 -7.7}

\arr{-4.7 -5.3}{-4.3 -5.7}
\arr{-2.7 -5.3}{-2.3 -5.7}
\arr{-.7 -5.3}{-.3 -5.7}
\arr{-5.7 -6.3}{-5.3 -6.7}
\arr{-3.7 -6.3}{-3.3 -6.7}
\arr{-1.7 -6.3}{-1.3 -6.7}
\arr{-6.7 -7.3}{-6.3 -7.7}
\arr{-4.7 -7.3}{-4.3 -7.7}
\arr{-2.7 -7.3}{-2.3 -7.7}
\arr{-.7 -7.3}{-.3 -7.7}

\arr{-8.3 -8.3}{-8.7 -8.7}
\arr{-7.7 -8.3}{-7.3 -8.7}
\arr{-6.3 -8.3}{-6.7 -8.7}
\arr{-5.7 -8.3}{-5.3 -8.7}
\arr{-4.3 -8.3}{-4.7 -8.7}
\arr{-3.7 -8.3}{-3.3 -8.7}
\arr{-2.3 -8.3}{-2.7 -8.7}
\arr{-1.7 -8.3}{-1.3 -8.7}
\arr{-.3 -8.3}{-.7 -8.7}

\arr{-9.3 -9.3}{-9.7 -9.7}
\arr{-8.7 -9.3}{-8.3 -9.7}
\arr{-7.3 -9.3}{-7.7 -9.7}
\arr{-6.7 -9.3}{-6.3 -9.7}
\arr{-5.3 -9.3}{-5.7 -9.7}
\arr{-4.7 -9.3}{-4.3 -9.7}
\arr{-3.3 -9.3}{-3.7 -9.7}
\arr{-2.7 -9.3}{-2.3 -9.7}
\arr{-1.3 -9.3}{-1.7 -9.7}
\arr{-.7 -9.3}{-.3 -9.7}

\setsolid

\linethickness.7mm
\putrule from 0 -3.3 to 0 -3.7
\putrule from 0 -4.3 to 0 -5.7
\putrule from 0 -6.3 to 0 -7.7
\putrule from 0 -8.3 to 0 -9.7

\putrule from -1 -3.3 to -1 -4.7
\putrule from -1 -5.3 to -1 -6.7
\putrule from -1 -7.3 to -1 -8.7
\putrule from -1 -9.3 to -1 -9.7

\put{$\ssize 3$} at 8 -4
\put{$\ssize 4$} at 8 -5
\put{$\ssize 5$} at 8 -6
\put{$\ssize 6$} at 8 -7
\put{$\ssize 7$} at 8 -8
\put{$\ssize 8$} at 8 -9

\put{$t$} at 8 -3

\setdots <.5mm>

\arr{.7 -5.3}{.3 -5.7}

\arr{1.7 -6.3}{1.3 -6.7}

\arr{2.7 -7.3}{2.3 -7.7}
\arr{.7 -7.3}{.3 -7.7}

\arr{.3 -4.3}{.7 -4.7}

\arr{1.3 -5.3}{1.7 -5.7}
\arr{2.3 -6.3}{2.7 -6.7}
\arr{.3 -6.3}{.7 -6.7}
\arr{3.3 -7.3}{3.7 -7.7}
\arr{1.3 -7.3}{1.7 -7.7}

\arr{4.3 -8.3}{4.7 -8.7}
\arr{3.7 -8.3}{3.3 -8.7}
\arr{2.3 -8.3}{2.7 -8.7}
\arr{1.7 -8.3}{1.3 -8.7}
\arr{.3 -8.3}{.7 -8.7}

\arr{5.3 -9.3}{5.7 -9.7}
\arr{4.7 -9.3}{4.3 -9.7}
\arr{3.3 -9.3}{3.7 -9.7}
\arr{2.7 -9.3}{2.3 -9.7}
\arr{1.3 -9.3}{1.7 -9.7}
\arr{.7 -9.3}{.3 -9.7}

\put{$1$} at 1 -5  
\put{$1$} at 2 -6  
\put{$1$} at 3 -7  
\put{$1$} at 4 -8  
\put{$1$} at 5 -9

\put{$6$} at 1 -7  
\put{$7$} at 2 -8
\put{$8$} at 3 -9

\put{$24$} at 0 -8  
\put{$32$} at 1 -9  

\comment
\setdots <1mm>
\plot -10.5 -5  -6 -5 /
\plot -10.5 -6  -7 -6 /
\plot -10.5 -7  -8 -7 /
\plot -10.5 -8  -9 -8 /

\put{$R_5$} at -10 -5.5
\put{$R_6$} at -10 -6.5
\put{$R_7$} at -10 -7.5
\endcomment

\multiput{$0$} at 3 -5  4 -6  5 -7  5 -5  6 -6  7 -7 /

\setsolid
\plot .6 -5.2  -1.6 -5.8 /
\plot 2.5 -5.1 -3.5 -5.9 /
\plot 4.5 -5.1 -5.5 -5.9 /

\plot -.2 -6.2  -.8 -6.8 /
\plot 1.5 -6.1 -2.5 -6.9 /
\plot 3.5 -6.1 -4.5 -6.9 /
\plot 5.5 -6.1 -6.5 -6.9 /

\plot .6 -7.2  -1.6 -7.8 /
\plot 2.5 -7.1 -3.5 -7.9 /
\plot 4.5 -7.1 -5.5 -7.9 /
\plot 6.5 -7.1 -7.5 -7.9 /


\put{$\ssize 1$} at -5  -6
\put{$\ssize 4$} at -3  -5.95
\put{$\ssize 7$} at -1.4  -5.9

\put{$\ssize 1$} at -6  -7
\put{$\ssize 5$} at -4  -7
\put{$\ssize 12$} at -2  -6.95
\put{$\ssize 12$} at -.4  -6.8

\put{$\ssize 1$} at -7  -8
\put{$\ssize 6$} at -5  -8
\put{$\ssize 18$} at -3  -7.95
\put{$\ssize 29$} at -1.4  -7.9

\endpicture}
$$
As we have mentioned already, of special interest seem to be the
knight's moves, which yield the numbers $1, 2, 7, 29, 130,\dots$
on the pylon of the even-index triangle.
	\bigskip
{\bf 4. Consequences.}
	\medskip
We see that {\it for the odd-index triangle, the entries on the left
determine those on the right, and conversely, the entries on the right determine
those on the left:}
	\medskip
{\bf Corollary 3.} 
$$
\align
 d''_i(t) &= d'_{i+2}(t+1)-2d'_{i+2}(t) + d'_{i+2}(t-1) \cr
 d'_i(t) &= d''_{i}(t+1)-2d''_{i-1}(t) + d''_{i-2}(t-1)
\endalign
$$
{\it The first equality holds for $i \le \frac{t-4}2$ (and $t \ge 4$), the second 
for $i < \frac t2$ (and $t\ge 1$),} thus for the vertices not lying on the pylons.
	\medskip
Proof: For the first equality, we use the formulae (2) and (3):
$$
\align
 d''_i(t) &= d_{i+1}(t)-d_{i+1}(t-1) \cr
          &=\bigl(d'_{i+2}(t+1)-d'_{i+2}(t)\bigr)-\bigl(d'_{i+2}(t)-d'_{i+2}(t-1)\bigr) \cr
          &= d'_{i+2}(t+1)-2d'_{i+2}(t)+d'_{i+2}(t-1).
\endalign
$$
Similarly, for the second equality, we use (1) and (4):
$$
\align
 d'_i(t) &= d_{i}(t)-d_{i-1}(t-1) \cr
          &= (d''_{i}(t+1)-d''_{i-1}(t))-(d''_{i-1}(t)-d''_{i-2}(t-1)) \cr
          &= d''_{i}(t+1)-2d''_{i-1}(t) + d''_{i-2}(t-1) 
\endalign
$$
	\medskip
These rules can be reformulated as follows, involving only
vertices in two consecutive layers:
$$
\align
 d''_i(t) &= 2d'_{i+1}(t) - d'_{i+1}(t-1) - d'_{i+2}(t)+d'_{i+2}(t-1)\tag{N} \cr
          &= d'_{i+2}(t+1) +  d'_{i+3}(t)- d'_{i+3}(t+1).\tag{N$'$}\cr \cr
 d'_i(t) &= 2d''_{i-2}(t) - d''_{i-1}(t+1) - d''_{i-1}(t) +d''_{i}(t-1)\tag{N} \cr
          &= d''_{i-2}(t-1) + d''_{i}(t) -  d''_{i-1}(t-1).\tag{N$'$}
\endalign
$$
Remark: Our interest in the vertices of two consecutive layers comes from the fact
that the dimension vectors of the modules $R_t$ are obtained by looking at two
consecutive layers.
	\medskip
Proof of the first N-rule: 
$$
\align
 d''_i(t) &= d'_{i+2}(t+1)-2d'_{i+2}(t) + d'_{i+2}(t-1) \cr
          &= \left(2'd_{i+1}(t) + d'_{i+2}(t)- d'_{i+1}(t-1)\right) 
     -2d'_{i+2}(t) + d'_{i+2}(t-1) \cr
          &= 2d'_{i+1}(t) - d'_{i+1}(t-1) -d'_{i+2}(t) + d'_{i+2}(t-1)
\endalign
$$
Proof of the first N$'$-rule:
$$
\align
 d''_i(t) &= d'_{i+2}(t+1)-2d'_{i+2}(t) + d'_{i+2}(t-1) \cr
          &= d'_{i+2}(t+1)-2d'_{i+2}(t) + \left(2d'_{i+2}(t) +                 d'_{i+3}(t)-d'_{i+3}(t+1)\right) \cr
          &= d'_{i+2}(t+1) +  d'_{i+3}(t)-d'_{i+3}(t+1)
\endalign
$$

The remaining rules are shown in the same way.
	\bigskip
There are also the following two summation formulae, adding up
the values along a sequence of south-east arrows starting at the left boundary.
In the case of the even-index triangle, the sequence has to stop before the
pylon. There is no such restriction in the case of the odd-index triangle.

	\medskip
{\bf Corollary 4. (Summation formulae)} 
	\smallskip

\item{(a)} {\it For all $0 \le i \le \frac{t-1}2$}
$$
 \sum_{0 \le j \le i}  d_j(t-i+j) = d''_i(t+1).
$$

\item{(b)} {\it For all $0 \le i$}
$$
 \sum_{0 \le j \le i} d'_j(t-i+j) = d_i(t).
$$
	\bigskip
In particular, since $d_i(t) = d_{t-i}(t)$, we see that
$$
 \sum_{0 \le j \le i} d'_j(t-i+j) = 
 \sum_{0 \le j \le t-i} d'_j(t-i+j).
$$
	\medskip

Proof: (a) If $i \le \frac{t-1}2,$ and $0\le j \le i$, then
$j \le \frac{t-1}2 \le \frac{t-i+j-1}2$, thus, according to (4),
$$
 d_j(t-i+j) = d_j''(t-i+j+1)-d''_{j-1}(t-i+j),
$$
and therefore
$$
\align
   \sum_{0 \le j \le i}  d_j(t-i+j) &= \sum_{0 \le j \le i}
   \left( d_j''(t-i+j+1)-d''_{j-1}(t-i+j)\right) \cr
  &=
   d_i''(t+1) - d''_{-1}(t-i) = d''_i(t+1).
\endalign
$$

Similarly, for (b) we use the formula (1). 
		\bigskip\bigskip
{\bf 5. Proof of Theorem.}
	\medskip
First, let us note that the rules (1) and (2) are actually equivalent,
due to the fact that $d_i(t) = d_{t-i}(t)$ and $d''_i(t) = d'_{t-i-1}(t)$.
Namely, assume (1) is satisfied. 
Then
$$
\align
 d''_i(t) &= d'_{t-i-1}(t)   \cr
          &= d_{t-i-1}(t)-d_{t-i-2}(t-1)  \cr
          &= d_{t-(t-i-1)}(t) - d_{(t-1)-(t-i-2)}(t-1) \cr
          &= d_{i+1}(t)-d_{i+1}(t-1),
\endalign
$$
thus (2) holds. In the same way, one sees the opposite implication.          
	\medskip

We are going to show that the theorem is 
a direct consequence of Proposition 4.1 in [FR2].
We fix some $t \ge 1$.
If $t$ is even, we assume that $x$ is a sink, otherwise that $x$ is 
a source. We choose the corresponding
bipartite orientation on the $3$-regular tree, this is the quiver whose
representations will be considered.
	\medskip
First, let us consider the rule (1).  Let $y$ be a neighbor of $x$.
According to [FR2], there is an exact
sequence 
$$
 0 \to P_{t-1}(y) \to P_t(x) \to R_t(x,y) \to 0,
$$
thus for any vertex $z$ of the quiver, there is an exact sequence
of vector spaces
$$
 0 \to P_{t-1}(y)_z \to P_t(x)_z \to R_t(x,y)_z \to 0.
$$
This means that
$$
 \dim R_t(x,y)_z = \dim P_t(x)_z - \dim P_{t-1}(y)_z.
$$
The only vertices $z$ to be considered are those with $D(x,z) \equiv t
\mod 2$.
We denote by $[x,z]$ the path between $x$ and $z$.
We have to distinguish whether $y$ belongs to $[x,z]$ or not.

First, consider the case $y\notin [x,z].$ Let $i = \frac12(t-D(x,z)),$
thus $2i \le t$ (since $t-2i = D(x,z) \ge 0$) and 
$$
   d'_i(t) = \dim R_t(x,y)_z = \dim P_t(x)_z - \dim P_{t-1}(y)_z.
$$
As we know, $\dim P_t(x)_z = d_i(t)$, thus it remains to calculate
$\dim P_{t-1}(y)_z$.
Since $y\notin [x,z]$, the path from $y$ to $z$ runs through $x$, thus
$$
 D(y,z) = 1+ D(x,z) = 1+t-2i = (t-1)-2(i-1),
$$
and consequently
$$
   \dim P_{t-1}(y)_z = d_{i-1}(t-1).
$$ 

Next, consider the case $y\in [x,z]$ and let $i = \frac12(D(x,z)+t),$ 
thus $D(x,z) = 2i-t$ and therefore
$$
 d'_i(t) = \dim R_t(x,y)_z.
$$
Since $y\in [x,z]$, we have 
$$
 D(y,z) = -1+D(x,z) = -1+2i-t,
$$
therefore
$$
 D(y,z) = |-1+2i-t| = |t-2i+1| = |(t-1)-2(i-1)|,
$$
but this means again that
$$
  \dim P_{t-1}(y)_z = d_{i-1}(t-1). 
$$

In both cases, we have shown that 
$$
 d'_i(t) = 
\dim R_t(x,y)_z = \dim P_t(x)_z - \dim P_{t-1}(y)_z = d_i(t)-d_{i-1}(t).
$$
This completes the proof of (1). 
	\medskip
Next we show the three rules (3), (4), (5). 

Let $y, y', y''$
be the neighbors of $x$. According to [FR2], there is an exact
sequence 
$$
 0 \to P_{t-1}(y') \to R_t(x,y) \to R_{t-1}(y'',x) \to 0,
$$
thus for any vertex $z$ of the quiver, there is an exact sequence
of vector spaces
$$
 0 \to P_{t-1}(y')_z \to R_t(x,y)_z \to R_{t-1}(y'',x)_z \to 0,
$$
This means that
$$
 \dim P_{t-1}(y')_z = \dim R_t(x,y)_z - \dim R_{t-1}(y'',x)_z. \tag{$*$}
$$
Again, we consider only the vertices $z$ with $D(x,z) \equiv t
\mod 2$.
	\medskip

For the proof of (3) and (4), we consider vertices $z$ with $y'\notin [x,z]$,
say with $D(x,z) = t-2i-2.$ 
We claim that 
$$
 \dim P_{t-1}(y')_z = d_i(t-1).
$$ 
In order to verify this equality, we have to show that
$D(y',z) = |t-1-2i|.$ 
But $D(y',z) = 1 + D(x,z) = 1+t-2(i+1) = t-1-2i.$
	\medskip
In order to establish the rule (3), we start with a vertex $z$
such that either $z = x$ or else $y''\in [x,z]$, thus in both cases
$y\notin [x,z]$. Let $i = \frac12(t-D(x,z))-1$, thus $D(x,z) = t-2i-2$
and 
$$
 d'_{i+1}(t) = \dim R_t(x,y)_z.
$$
Next, let us show that 
$$
 d'_{i+1}(t-1) = \dim R_{t-1}(y'',x)_z
$$
This holds true in case $z = x$: namely, then $x\in [y'',z]$ and
$D(y'',x) = 1 = 2(i+1)-(t-1)$ (since $0 = D(x,z) = t-2(i+1)$).
If $z \neq x$, then 
$x\notin [y'',z]$ and $D(y'',z) = D(x,z) -1 = t-2(i+1)-1 = 
(t-1) - 2(i+1)$.

Altogether, we see that $(*)$ yields the required equality
$$  
 d_i(t-1) =   d'_{i+1}(t) -  d'_{i+1}(t-1).
$$

Let us discuss the rule (4). 
The case $i = \frac{t-2}2$ follows from (3), since for this value of $i$,
we have both $d''_i(t) = d'_{i+1}(t)$ as well as 
$d''_{i-1}(t-1) = d'_{i+1}(t-1)$. Thus, we only have to consider 
the cases $i < \frac{t-2}2.$

We start with a vertex $z$ such that $y\in [x,z].$ Let 
$i = \frac12(t-D(x,z)-2)$, thus $D(x,z) = t-2i-2$ and
$$
 d''_i(t) = \dim R_t(x,y)_z.
$$
Note that $x \in [y'',z]$ and $D(y'',z) = 1+D(x,z) = t-2i-1 = (t-1)-2(i-1)-2$,
thus
$$
 d''_{i-1}(t-1) = \dim R_{t-1}(y'',x)_z.
$$
It remains to show that $\dim P_{t-1}(y')_z = d_i(t-1).$

Since $D(x,z) = t-2i-2$, we know that $\dim P_{t-1}(y')_z = d_i(t-1),$ 
thus $(*)$ yields the required equality
$$  
 d_i(t-1) =   d''_{i}(t) -  d''_{i-1}(t-1).
$$

The final considerations concern the rule (5). The case $i = \frac t2$ of (5)
coincides with the case $j = \frac{t-2}2$ of (3) and (4), thus we only have
to establish the cases $i < \frac t2.$ 

We consider a vertex $z$ such that
$y' \in [x,z]$ and let $i = \frac12(t-D(x,z))$, thus $D(x,z) = t-2i$
and $D(y',z) = -1+D(x,z) = t-2i-1.$ 
Then 
$$
 \dim P_{t-1}(y')_z = d_i(t-1),
$$ 
since $D(y',z) = t-1-2i.$
We have
$$
 d'_i(t) = \dim R_t(x,y)_z,
$$
since $y\notin [x,z]$ and $D(x,z) = t-2i$.
And we have
$$
 d''_{i-2}(t-1) = \dim R_{t-1}(y'',x)_z,
$$
since $x\in [y'',z]$ and $D(y'',z) = 1+D(x,z) = t-2i+1 = 
(t-1)-2(i-2)-2.$
Altogether, we see again that $(*)$ yields the required equality, namely now
$$  
 d_i(t-1) =   d'_{i}(t) -  d''_{i-2}(t-1).
$$
This completes the proof.
	\bigskip\bigskip
\vfill\eject
{\bf 6. Questions and remarks.}
	\medskip
{\bf 6.1. The south-west arrows of the odd-index triangle.} 
The relationship between the two triangles
was established using differences between numbers, mostly along arrows.
In this way, all the arrows of the even-index triangle, and all
the south-west arrows of the odd-index triangle have been used.

It seems to be of interest to understand also the meaning of the differences
along the sout-east arrows in the odd-index triangle. Here are the numbers 
$d'_i(t)-d'_{i-1}(t-1)$, for $0\le t \le 12$:
	\medskip
$$
{\beginpicture
\setcoordinatesystem units <.5cm,.7cm>
\multiput{$1$} at 1 1  2 2 
  0 0  1 -1   3 -3  4 -4  5 -5  6 -6  7 -7  8 -8  9 -9  
   10 -10   -1 -1  -2 -2  -3 -3 -4 -4  -5 -5  -6 -6  -7 -7 -8 -8  -9 -9  
   -10 -10  /
\put{$2$} at 2 -2
\put{$2$} at 0 -2
\put{$3$} at -1 -3
\put{$5$} at 1 -3
\put{$4$} at -2 -4
\put{$9$} at 0 -4
\put{$9$} at 2 -4
\put{$5$} at -3 -5
\put{$14$} at -1 -5
\put{$22$} at 1 -5
\put{$7$} at 3 -5
\put{$6$} at -4 -6
\put{$20$} at -2 -6
\put{$41$} at 0 -6
\put{$42$} at 2 -6
\put{$8$} at  4 -6
\put{$7$} at -5 -7
\put{$27$} at -3 -7
\put{$67$} at -1 -7
\put{$102$} at 1 -7
\put{$40$} at 3 -7
\put{$9$} at 5 -7
\put{$8$} at -6 -8
\put{$35$} at -4 -8
\put{$101$} at -2 -8
\put{$195$} at 0 -8
\put{$202$} at 2  -8
\put{$50$} at 4 -8
\put{$10$} at 6 -8
\put{$9$} at -7 -9
\put{$44$} at -5 -9
\put{$144$} at -3 -9
\put{$330$} at -1 -9
\put{$490$} at 1 -9
\put{$217$} at 3 -9
\put{$61$} at 5 -9
\put{$11$} at 7 -9
\put{$10$} at -8 -10
\put{$54$} at -6 -10
\put{$197$} at -4 -10
\put{$517$} at -2 -10
\put{$955$} at 0 -10
\put{$995$} at 2 -10
\put{$289$} at 4 -10
\put{$73$} at 6 -10
\put{$12$} at 8 -10

\endpicture}
$$
	\medskip
{\bf 6.2. The difference operator for sequences.}
Some of the results presented here, may be reformulated 
using the following operators $\Delta$ and $E$.
Given a function $u\:\Bbb N_0 \to \Bbb R$, we define new functions on $\Bbb N_0$
as follows
$$
\align
 (\Delta u)(t) &= u(t+1)-u(t) \cr
 (Eu)(t) &= u(t+1),
\endalign
$$
the operator $\Delta$ is called the difference operator, $E$ the shift operator 
(see [E]).
For example, the sequence $f$ of Fibonacci numbers satisfies
the condition 
$\Delta E f = f.$ Let us stress that we have used the difference operator
in the proof of the Corollaries. Here we note that the 
assertions (2) and (3) concern the difference and the shift operators, as follows:
$$
\align
 E d''_i &= \Delta d_{i+1}, \cr
 d_i &= \Delta d'_{i+1}. 
\endalign
$$
	
	\medskip
The first assertion of Corollary 3 concerns the following identity:
$$
 Ed''_i = \Delta d_{i+1} = \Delta^2 d'_{i+2}.
$$
	\bigskip
{\bf 6.3. Group actions on quivers and valued quivers.}
	\medskip
The valued quivers which have been considered in the paper are
derived from group actions an the $3$-regular tree $T$. If $x,y$ is
a pair of neighboring vertices, we may consider the following
groups of automorphisms of $T$. Let $G_x$ denote the group of
automorphisms which fix $x$, let $G_{xy}$ denote the group
of automorphisms which fix both $x$ and $y$. Then the functions on
$T_0$ which are $G_x$-invariant may be identified with the functions
on $\Delta^{\ev},$  whereas the functions on
$T_0$ which are $G_{xy}$-invariant may be identified with the functions
on $\Delta^{\odd}.$ 
	\bigskip
{\bf 6.4. Left hammocks.}
	\medskip
We should mention that the additive functions $d, d'$ 
exhibited here are typical left hammock functions, see [RV].
	\bigskip
{\bf 6.5. Delannoy paths.}
	\medskip
Hirschhorn [H] has shown that the second last 
column of the even-index triangle (with numbers $1,3,12,53,247,1192,\dots$) is just 
the sequence A110122 in Sloane's On-Line Encyclopedia of Integer Sequences, it
counts the number of the Delannoy paths from $(0,0)$ to $(n,n)$ which do not cross
horizontally the diagonal $x=y$ (we recall that a Delannoy path is a sequence of
steps $(1,0), (1,1), (0,1)$ in the plane (thus going north, northeast and east);
and such a path is said to cross the diagonal $x=y$ horizontally provided it contains
a subpath of the form $(m-1,m) \to (m,m) \to (m+1,m)$). 

Hirschhorn's proof is computational and does not provide an intrinsic
relationship between the Delannoy paths in question and say suitable 
elements of the
Fibonacci modules. It seems to be of interest to establish a direct
relationship.
	\bigskip
	\vfill\eject
{\bf 6.6. Concordance.}
	\medskip
As we have mentioned, the numbers displayed in the even-index triangle have been
considered already in the paper [FR1],
as well in subsequent publications by other authors. In [FR1], these numbers 
have also been
denoted by $b_i[j]$ and $c_i[j]$; for the convenience of the reader, let us
present the different notations for the numbers of the triangle:
$$
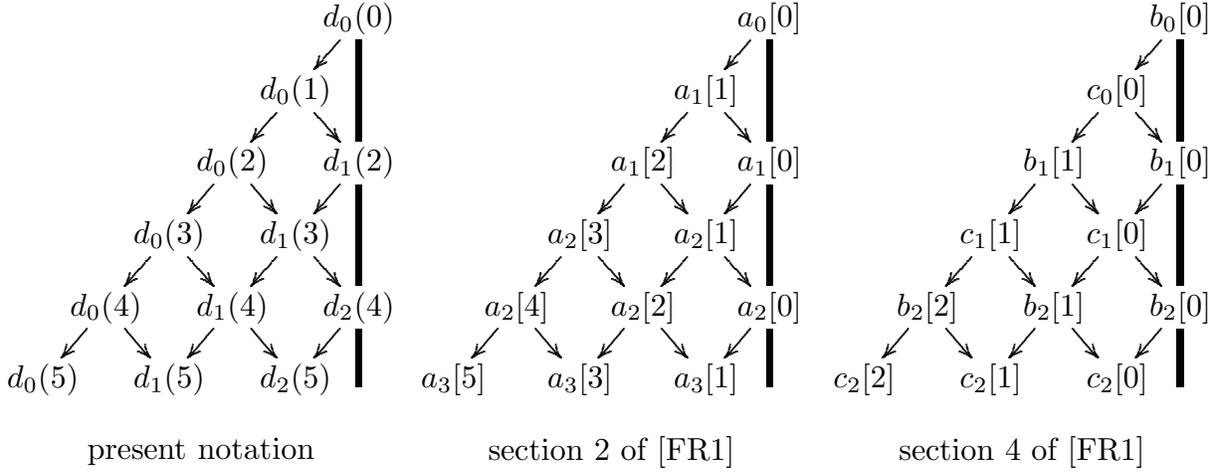
{\beginpicture
\setcoordinatesystem units <.7cm,.8cm>
\put{\beginpicture
\put{$d_0(0)$} at 0 0 
\put{$d_0(1)$} at -1 -1
\put{$d_0(2)$} at -2 -2
\put{$d_1(2)$} at 0 -2
\put{$d_0(3)$} at -3 -3
\put{$d_1(3)$} at -1 -3
\put{$d_0(4)$} at -4 -4
\put{$d_1(4)$} at -2 -4
\put{$d_2(4)$} at 0 -4 
\put{$d_0(5)$} at -5 -5
\put{$d_1(5)$} at -3 -5
\put{$d_2(5)$} at -1 -5 
\arr{-.3 -.3}{-.7 -.7}

\arr{-1.3 -1.3}{-1.7 -1.7}
\arr{-.7 -1.3}{-.3 -1.7}

\arr{-2.3 -2.3}{-2.7 -2.7}
\arr{-1.7 -2.3}{-1.3 -2.7}
\arr{-.3 -2.3}{-.7 -2.7}

\arr{-3.3 -3.3}{-3.7 -3.7}
\arr{-2.7 -3.3}{-2.3 -3.7}
\arr{-1.3 -3.3}{-1.7 -3.7}
\arr{-.7 -3.3}{-.3 -3.7}

\arr{-4.3 -4.3}{-4.7 -4.7}
\arr{-3.7 -4.3}{-3.3 -4.7}
\arr{-2.3 -4.3}{-2.7 -4.7}
\arr{-1.7 -4.3}{-1.3 -4.7}
\arr{-.3 -4.3}{-.7 -4.7}

\put{} at 0 -5
\setsolid

\linethickness.7mm
\putrule from 0 -.3 to 0 -1.7
\putrule from 0 -2.3 to 0 -3.7
\putrule from 0 -4.3 to 0 -5.1
\endpicture} at 0 0
\put{\beginpicture
\put{$a_0[0]$} at 0 0 
\put{$a_1[1]$} at -1 -1
\put{$a_1[2]$} at -2 -2
\put{$a_1[0]$} at 0 -2
\put{$a_2[3]$} at -3 -3
\put{$a_2[1]$} at -1 -3
\put{$a_2[4]$} at -4 -4
\put{$a_2[2]$} at -2 -4
\put{$a_2[0]$} at 0 -4 
\put{$a_3[5]$} at -5 -5
\put{$a_3[3]$} at -3 -5
\put{$a_3[1]$} at -1 -5 
\arr{-.3 -.3}{-.7 -.7}

\arr{-1.3 -1.3}{-1.7 -1.7}
\arr{-.7 -1.3}{-.3 -1.7}

\arr{-2.3 -2.3}{-2.7 -2.7}
\arr{-1.7 -2.3}{-1.3 -2.7}
\arr{-.3 -2.3}{-.7 -2.7}

\arr{-3.3 -3.3}{-3.7 -3.7}
\arr{-2.7 -3.3}{-2.3 -3.7}
\arr{-1.3 -3.3}{-1.7 -3.7}
\arr{-.7 -3.3}{-.3 -3.7}

\arr{-4.3 -4.3}{-4.7 -4.7}
\arr{-3.7 -4.3}{-3.3 -4.7}
\arr{-2.3 -4.3}{-2.7 -4.7}
\arr{-1.7 -4.3}{-1.3 -4.7}
\arr{-.3 -4.3}{-.7 -4.7}

\put{} at 0 -5
\setsolid

\linethickness.7mm
\putrule from 0 -.3 to 0 -1.7
\putrule from 0 -2.3 to 0 -3.7
\putrule from 0 -4.3 to 0 -5.1
\endpicture} at 6.5 0
\put{\beginpicture
\put{$b_0[0]$} at 0 0 
\put{$c_0[0]$} at -1 -1
\put{$b_1[1]$} at -2 -2
\put{$b_1[0]$} at 0 -2
\put{$c_1[1]$} at -3 -3
\put{$c_1[0]$} at -1 -3
\put{$b_2[2]$} at -4 -4
\put{$b_2[1]$} at -2 -4
\put{$b_2[0]$} at 0 -4 
\put{$c_2[2]$} at -5 -5
\put{$c_2[1]$} at -3 -5
\put{$c_2[0]$} at -1 -5 
\arr{-.3 -.3}{-.7 -.7}

\arr{-1.3 -1.3}{-1.7 -1.7}
\arr{-.7 -1.3}{-.3 -1.7}

\arr{-2.3 -2.3}{-2.7 -2.7}
\arr{-1.7 -2.3}{-1.3 -2.7}
\arr{-.3 -2.3}{-.7 -2.7}

\arr{-3.3 -3.3}{-3.7 -3.7}
\arr{-2.7 -3.3}{-2.3 -3.7}
\arr{-1.3 -3.3}{-1.7 -3.7}
\arr{-.7 -3.3}{-.3 -3.7}

\arr{-4.3 -4.3}{-4.7 -4.7}
\arr{-3.7 -4.3}{-3.3 -4.7}
\arr{-2.3 -4.3}{-2.7 -4.7}
\arr{-1.7 -4.3}{-1.3 -4.7}
\arr{-.3 -4.3}{-.7 -4.7}

\put{} at 0 -5
\setsolid

\linethickness.7mm
\putrule from 0 -.3 to 0 -1.7
\putrule from 0 -2.3 to 0 -3.7
\putrule from 0 -4.3 to 0 -5.1
\endpicture} at 13 0

\put{present notation} at 0  -3.5
\put{section 2 of [FR1]} at 6.5  -3.5
\put{section 4 of [FR1]} at 13  -3.5

\endpicture}
$$

Similarly, let us consider the odd-index triangle. Again, 
for the convenience of the reader, let us compare the
new notation (left) with the notation used in the paper [FR2] (right):
$$
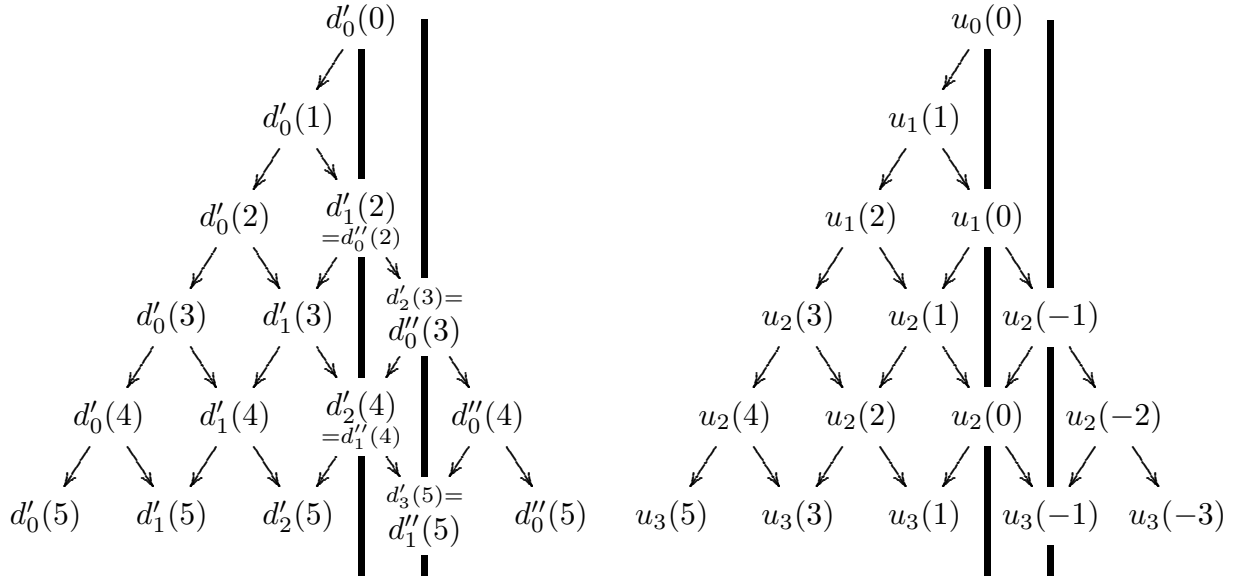
{\beginpicture
\setcoordinatesystem units <.7cm,1.1cm>

\put{\beginpicture
\put{$d'_0(0)$} at -1 -1 
\put{$d'_0(1)$} at -2 -2 
\put{$d'_0(2)$} at -3 -3 
\put{$d'_0(3)$} at -4 -4 
\put{$d'_0(4)$} at -5 -5 
\put{$d'_0(5)$} at -6 -6

\put{$d'_1(2)$} at -1 -2.9 
\put{$\ssize = d''_0(2)$} at -1 -3.2
\put{$d'_1(3)$} at -2 -4 
\put{$d'_1(4)$} at -3 -5 
\put{$d'_1(5)$} at -4 -6 

\put{$\ssize d'_2(3) = $} at 0 -3.8
\put{$d''_0(3)$} at 0 -4.15
\put{$d'_2(4)$} at -1 -4.9
\put{$\ssize = d''_1(4)$} at -1 -5.2
\put{$d'_2(5)$} at -2 -6 

\put{$\ssize d'_3(5)= $} at 0 -5.8
\put{$d''_1(5)$} at 0 -6.15

\arr{-1.3 -1.3}{-1.7 -1.7}

\arr{-2.3 -2.3}{-2.7 -2.7}
\arr{-1.7 -2.3}{-1.4 -2.6}

\arr{-3.3 -3.3}{-3.7 -3.7}
\arr{-2.7 -3.3}{-2.3 -3.7}
\arr{-1.4 -3.4}{-1.7 -3.7}
\arr{-.6 -3.4}{-.4 -3.6}

\arr{-4.3 -4.3}{-4.7 -4.7}
\arr{-3.7 -4.3}{-3.3 -4.7}
\arr{-2.3 -4.3}{-2.7 -4.7}
\arr{-1.7 -4.3}{-1.4 -4.6}
\arr{-.4 -4.4}{-.6 -4.6}

\arr{-5.3 -5.3}{-5.7 -5.7}
\arr{-4.7 -5.3}{-4.3 -5.7}
\arr{-3.3 -5.3}{-3.7 -5.7}
\arr{-2.7 -5.3}{-2.3 -5.7}
\arr{-1.4 -5.4}{-1.7 -5.7}
\arr{-.6 -5.4}{-.4 -5.6}

\linethickness.7mm
\putrule from 0 -1 to 0 -3.6
\putrule from 0 -4.4 to 0 -5.6
\putrule from 0 -6.4 to 0 -6.6

\putrule from -1 -1.3 to -1 -2.6
\putrule from -1 -3.4 to -1 -4.6
\putrule from -1 -5.4 to -1 -6.6

\arr{.4 -4.4}{.7 -4.7}

\arr{1.3 -5.3}{1.7 -5.7}
\arr{.7 -5.3}{.4 -5.6}

\put{$d''_0(4)$} at 1 -5  
\put{$d''_0(5)$} at 2 -6  

\endpicture} at 0 0
\put{\beginpicture
\put{$u_0(0)$} at -1 -1 
\put{$u_1(1)$} at -2 -2 
\put{$u_1(2)$} at -3 -3 
\put{$u_2(3)$} at -4 -4 
\put{$u_2(4)$} at -5 -5 
\put{$u_3(5)$} at -6 -6

\put{$u_1(0)$} at -1 -3 
\put{$u_2(1)$} at -2 -4 
\put{$u_2(2)$} at -3 -5 
\put{$u_3(3)$} at -4 -6 

\put{$u_2(-1)$} at 0 -4
\put{$u_2(0)$} at -1 -5
\put{$u_3(1)$} at -2 -6 

\put{$u_3(-1)$} at 0 -6

\arr{-1.3 -1.3}{-1.7 -1.7}

\arr{-2.3 -2.3}{-2.7 -2.7}
\arr{-1.7 -2.3}{-1.3 -2.7}

\arr{-3.3 -3.3}{-3.7 -3.7}
\arr{-2.7 -3.3}{-2.3 -3.7}
\arr{-1.3 -3.3}{-1.7 -3.7}
\arr{-.7 -3.3}{-.3 -3.7}

\arr{-4.3 -4.3}{-4.7 -4.7}
\arr{-3.7 -4.3}{-3.3 -4.7}
\arr{-2.3 -4.3}{-2.7 -4.7}
\arr{-1.7 -4.3}{-1.3 -4.7}
\arr{-.3 -4.3}{-.7 -4.7}

\arr{-5.3 -5.3}{-5.7 -5.7}
\arr{-4.7 -5.3}{-4.3 -5.7}
\arr{-3.3 -5.3}{-3.7 -5.7}
\arr{-2.7 -5.3}{-2.3 -5.7}
\arr{-1.3 -5.3}{-1.7 -5.7}
\arr{-.7 -5.3}{-.3 -5.7}

\linethickness.7mm
\putrule from 0 -1 to 0 -3.7
\putrule from 0 -4.3 to 0 -5.7
\putrule from 0 -6.3 to 0 -6.6

\putrule from -1 -1.3 to -1 -2.7
\putrule from -1 -3.3 to -1 -4.7
\putrule from -1 -5.3 to -1 -6.6

\arr{.3 -4.3}{.7 -4.7}

\arr{1.3 -5.3}{1.7 -5.7}
\arr{.7 -5.3}{.3 -5.7}

\put{$u_2(-2)$} at 1 -5  
\put{$u_3(-3)$} at 2 -6  

\endpicture} at 10 0
\endpicture}
$$
	\bigskip\bigskip

{\bf 7. References}
	\medskip
\item{[E]} S. N. Elaydi: An Introduction to Difference Equations.
  Springer (1996).
\item{[F]} Ph. Fahr.
   Infinite Gabriel-Roiter measures for the 3-Kronecker quiver.
  Dissertation Bielefeld. (2008).
\item{[FR1]} Ph. Fahr, C.M. Ringel: A partition formula for Fibonacci numbers, J.
   Integer Sequences, 11 (2008), Paper 08.1.4.
\item{[FR2]} Ph. Fahr, C.M. Ringel: Categorification of the Fibonacci Numbers
   Using Representations of Quivers. Preprint 2011.
\item{[H]} M.D. Hirschhorn: On Recurrences of Fahr and Ringel Arising in Graph
  Theory. J. Integer Sequences, 12 (2009), Paper 09.6.8.
\item{[HPR]} D. Happel, U. Preiser, C.M. Ringel: Vinberg's characterization of Dynkin
  diagrams using subadditive functions with applications to DTr-periodic modules.
  In: Representation Theory II (ed V. Dlab, P. Gabriel), 
  Springer LNM 832 (1980),     280-294.
 \item{[RV]} C.M. Ringel, D. Vossieck: Hammocks. Proc. London Math. Soc. (3) 54   (1987), 216-246.

	\bigskip\bigskip

{\rmk Philipp Fahr\par
Brehmstr. 51\par
D-40239 D\"usseldorf, Germany\par
e-mail: {\ttk philfahr\@gmx.net}\par
	\medskip
Claus Michael Ringel \par
Fakult\"at f\"ur Mathematik, Universit\"at Bielefeld \par
POBox 100\,131, \ D-33\,501 Bielefeld, Germany \par
e-mail: {\ttk ringel\@math.uni-bielefeld.de} \par}

\bye